\documentclass[11pt]{article}
\usepackage{amssymb}
\usepackage{amsfonts}
\usepackage{amscd}

\usepackage[english]{babel}

\usepackage{amsmath,latexsym,enumerate,amsthm,amssymb}

\newcommand{\be}{\begin{enumerate}}
\newcommand{\ee}{\end{enumerate}}
\let\ds=\displaystyle

\def\N{{\mathbb N}} \def\Z{{\mathbb Z}}

 \def\Q{{\mathbb Q}}

\def\R{{\mathbb R}} \def\C{{\mathbb C}}

\def\Sp{{\mathbb S}}

\def\T{{\cal T}}

\def\s{{\bf s}}

\def\m{{\bf m}}
\def\i{{\bf i}}

\def\z{{\bf z}}
\def\e{{\bf e}}
\def\x{{\bf x}}
\def\y{{\bf y}}
\def\a{{\bf a}}
\def\cb{{\bf c}}

\def\u{{\bf u}}
\def\v{{\bf v}}

\def\b{{\bf b}}
\def\m{{\bf m}}
\def\X{{\bf X}}

\def\h{{\bf h}}
\def\f{{\bf f}}

\newcommand{\zerob}{\boldsymbol{0}}
\newcommand{\unb}{\boldsymbol{1}}

\newcommand{\alphab}{\boldsymbol{\alpha}}
\newcommand{\betab}{\boldsymbol{\beta}}
\newcommand{\gammab}{\boldsymbol{\gamma}}

\newcommand{\mub}{\boldsymbol{\mu}}
\newcommand{\nub}{\boldsymbol{\nu}}
\newcommand{\sigmab}{\boldsymbol{\sigma}}
\newcommand{\taub}{\boldsymbol{\tau}}
\newcommand{\lambdab}{\boldsymbol{\lambda}}
\newcommand{\rhob}{\boldsymbol{\rho}}
\newcommand{\etab}{\boldsymbol{\eta}}

\def\cqfd{ $\diamondsuit $ }

\def\eps{{ \varepsilon }}

\catcode`\ç=13
\defç{\c{c}}
\catcode`\é=13
\defé{\'e}
\catcode`\à=13
\defà{\`a}
\catcode`\è=13
\defè{\`e}
\catcode`\â=13
\defâ{\^a}
\catcode`\ù=13
\defù{\`u}
\catcode`\ê=13
\defê{\^e}
\catcode`\î=13
\defî{\^\i}
\catcode`\ô=13
\defô{\^o}

\newcommand{\la}{\langle}

\newcommand{\ra}{\rangle}

\newtheorem{Lemme}{Lemma}

\newtheorem{Corollaire}{Corollary}

\newtheorem{Proposition}{Proposition}

\newtheorem{Theoreme}{Theorem}

\newtheorem{Definition}{Definition}

\begin{document}
\title{Mixed zeta functions and application to some
  lattice points problems.} 
\author{D. Essouabri\footnote{ 
Universit\'e de Caen,
 UFR des Sciences, Campus 2,
 Laboratoire de Math\'ematiques Nicolas Oresme (CNRS UMR 6139),
 Bd. Mar\'echal Juin, B.P. 5186,
14032 Caen, France.
Email : essoua@math.unicaen.fr}}

\date{ }
\maketitle 

\noindent
{\bf {\large Abstract:}}
{\it We consider zeta functions: $Z(f ;P ;s)=\sum_{\m \in \N^{n}} f(m_1,\dots,
  m_n)~P(m_1,\dots ,m_n)^{-s/d}$ where $P \in \R [X_1,\dots ,X_n]$ has degree $d$
  and $f$ is a function arithmetic in origin, e.g. a multiplicative
function. In this paper, I study the 
 meromorphic continuation of such series beyond an a priori domain of absolute convergence when   $f$ and $P$ satisfy 
 properties one typically meets in applications.  As a result, I prove an explicit asymptotic  for a general  class of  lattice
point problems  subject to  arithmetic   constraints.}\par  

\vskip 0.5cm

{\bf Mathematics Subject Classifications: 11M41, 11P21, 11N25, 11N37.}\par
{\bf Key words: Zeta functions, meromorphic continuation, Newton
  polyhedron, multiplicative functions, Lattice points, representation of integers.}
\vskip .2 in
\setcounter{tocdepth}{2}
\tableofcontents
 
\section{Introduction} 
Let $f:\N^{n} \rightarrow \C$ be a function and $P\in \R[X_1,\dots,X_n]$
a polynomial of degree $d$.
We define a mixed zeta function (associated to the pair $(f, P)$) as the series defined formally by: 
$$s\mapsto Z(f;P;s)=\sum_{\m \in \N^{n}} \frac{f(m_1,\dots,
  m_n)}{P(m_1,\dots ,m_n)^{s/d}} \quad (s \in \C).$$
These generating functions  are natural objects to study  when it is of interest 
to understand the asymptotic   density   of $f$ (on average) restricted to   the increasing family of
sets $\{P(\m)^{1/d} \le t\}$ as $t\to \infty.$ 
To do this, a classical method tells us that we must first find the domain of convergence $D$ of $Z(f;P;s)$
and then  understand its behavior along the boundary of $D.$\par 
In the classical case when $f$ is a polynomial
(eventually twisted by additive characters), after work of many authors, the problem is now 
understood for a large class of   $P$ (see \cite{mellin},
\cite{mahler}, \cite{pierrette}, \cite{sargosfourier},
\cite{sargosthese}, \cite{lichtincompo}, \cite{essouabrifourier}).
Some results have also been obtained 
if $f$ is the characteristic function of a suitable open
semi-algebraic subset of $\R^n$ (see \cite{mahlerquadratique}, 
\cite{essouabricompo}).\par
However, for many standard arithmetic problems, the function $f$ is  irregular; typically, though not always,  it will be
multiplicative. In such cases, no  general methods
to study $Z(f;P;s)$ are known  since the works cited above need to begin with 
  an integral representation for the series. This requires 
  the function  $f$ to have some reasonable expression as an algebraic or analytic function. If $f$ is
multiplicative, then it can sometimes occur that   methods that begin with an Euler product expression for the series (see
for example \cite{bretechecompo} and \cite{eulerprod}) can  be used,
but in such cases, $P$ must  also be multiplicative, that is, a monomial.\par 
The  point of this article is that it gives a   method to study the analytic properties for  
``mixed" zeta functions $Z(f;P;s),$ where $f$ can be an irregular function, and $P$ is any
homogeneous polynomial with positive coefficients (which should suffice for typical
arithmetical applications). We call the  class of   $f$  ``functions
 of  finite type" (see
\S 3.1). Such functions are often encountered in a variety of arithmetic problems.  For applications, it is quite useful to
have as precise information as possible about a first pole of $Z(f;P;s).$ Our main result, Theorem \ref{partieprincipale} (\S 3.2),  gives  a
  criterion, which if satisfied,   explicitly identifies the first pole as well as the leading term of the  principal
part  of $Z(f;P;s)$ at that pole.\par 
The applications  that we present to illustrate this theorem are stated 
in \S 3.3, and proved in \S 4.3ff. They illustrate, in particular, the type of  result that can be proved  whenever   the 
characteristic function $1_B$ of a subset of  
$\N^{n}$ is known (or can be shown) to be of finite type.  We give one
simple example here that is a straightforward consequence of Corollary
\ref{applicationneuve}.  Let \\
$\ds \quad B = \{\m = (m_1,\dots, m_n) \in \N^n : \ \forall i ~m_i
{\mbox { is square free }} \}. $ \\ 
When $n = 1,$   it is a classical fact (see
\cite{ten}, Chap.  I.3.7) that the number of squarefree integers in the interval $[1, t]$ is  asymptotic to $(6/\pi^2)
t.$ Thus,   $\# B \cap [1, t]^n$ is known to be asymptotic to $(6/\pi^2)^n t^n.$ If  we now change the enclosing region
from a box to one that is curved, say  $\{P^{1/d} \le t\} \cap \R_+^n,$ we would expect $\# B \cap \{P^{1/d} \le t\}$ to
grow at a rate
$C_n(B, P) t^n\,,$ where the  constant reflects both   arithmetic and geometric properties of this set. Indeed, our result
shows that {\it for any elliptic
polynomial $P$ of degree $d$}, the two features are independent of one another in the following sense.\\ 
 There exists 
$\theta >0$ such that   as  $t\rightarrow +\infty:$ 
 $ \ds \# B \cap \{ P^{1/d}(\m)\leq t\}= C_n(B, P)~ { t^{n}} +O\big(t^{ {n} -\theta}\big),$ \\ 
where $C_n(B, P):= \left(\frac{6}{\pi^2}\right)^n \cdot \left(\frac{1}{n} 
~\int_{\Sp^{n-1} \cap \R_+^n} P_d^{-n/d}(\v)~ d\sigma(\v) \right)$ (
$\Sp^{n-1}$ is the unit sphere).\\ 
The factor in parentheses is geometric and first appeared in the
lattice point problem studied by Mahler \cite{mahler}.\par  
The starting point of our method is the classical and remarkable
following Mellin's formula: 
\begin{equation}\label{mellinformula}
\frac{\Gamma (s)}{\left(\sum_{k=0}^r w_k \right)^{s}}
= 
\frac{1}{(2\pi i)^r} 
\int_{\rho_1-i\infty}^{\rho_1+i\infty}\dots
\int_{\rho_r-i\infty}^{\rho_r+i\infty}  
\frac{\Gamma(s -z_1-\dots-z_r)~\prod_{i=1}^r \Gamma (z_i) ~d\z}
{w_0^{s-z_1-\dots-z_r }\left(\prod_{k=1}^r w_k^{z_k}\right)}
\end{equation}
valid if $\forall i=0,\dots, r$, $\Re (w_i) >0$, 
$\forall i=1,\dots, r$ $\rho_i >0$ and $\Re(s) >\rho_1+\dots +\rho_r$.\par
This formula implies that for $\rhob \in \R_+^{*n}$ and  
$\Re(s) \gg 1$, 
$$Z(f;P;s)= 
\frac{1}{(2\pi i)^r} 
\int_{\rho_1-i\infty}^{\rho_1+i\infty}\dots
\int_{\rho_r-i\infty}^{\rho_r+i\infty}  F(s;z_1,\dots,z_r) ~dz_1\dots
dz_r,$$
where $F$ is a meromorphic function in an open subset of $\C^{r+1}$ 
(see \S 4.1 for more details).\par    
In the classical case (i.e. $f$ is a polynomial possibly twisted by additive
characters)   Mellin's formula was used by many authors
(\cite{mellin}, \cite{pierrette},\cite{matsumotosurvey1},...). 
In this event,  the function $F$ {\it has a meromorphic continuation to the
  whole space $\C^{r+1}$}. Applying   induction  
on the number of monomials of $P$, one then concludes that $Z(f;P;s)$ has a meromorphic
continuation to $\C$. This method gives, except for some very
special cases, only a set of possible poles. So, one cannot yet  use it to determine the dominant term in the principal
part at the first pole.\par 
If $f$ is of finite type, the function that plays the role of $F$ will {\it not   have}, in general, a
meromorphic continuation to  $\C^{r+1}$. Indeed, for many multiplicative functions, 
this is known to be impossible!  Consequently the method that works in
the classical case cannot be immediately applied to study
$Z(f;P;s)$ outside an a priori domain of absolute convergence, given,
say, by $\Re s > c$. In particular, one can
only expect to obtain a meromorphic continuation  of $Z(f;P;s)$ to a half-plane of the
form $\{\Re(s)>\eta\}$ along whose boundary there are essential
singularities. Our method identifies a possible first
pole $\sigma_0 \le c,$ in complete generality, from which an upper 
bound for the counting function 
$N(f;P;t):=\sum_{\{\m \in \N^{n}; ~P(\m)^{1/d}\leq t\}}
f(m_1,\dots,m_n)$ immediately follows. \par 
The particular interest of Theorem \ref{partieprincipale} is its proof
that the candidate first   pole is a
genuine pole of a precisely given order, provided that a certain
analytic criterion is satisfied, from which, of course,
follows the  explicit  dominant term for the asymptotic of $N(f;P;t).$  
Verification of this criterion requires some additional information
about the behavior of the multivariable Dirichlet series
${\cal M}(f;\s)$ (see \S 3.1) in small neighborhoods of certain points 
on the boundary of its domain of analyticity.
Such information is quite similar to that which was needed in  earlier
work  on multivariable Tauberian theorems by Lichtin \cite{lichtinduke1}
and de la Bretèche \cite{bretechecompo}. This criterion can be verified in specific
cases, as our examples indicate.\par 
The principal idea in the proof of Theorem 1 (see \S 4.2) is to 
associate in a natural way  several
(``mixed'') invariants to $f$ and $P$. These combinatorial-geometric
invariants allow  a good control of the data within an induction argument, and, in
particular,  play a very important role in the proofs of the crucial 
lemmas \ref{fonda1} and \ref{fonda3} (see \S 4.1).
\section{Preliminaries}
\subsection{Notations}
\be
\item 
$\N=\{1,2,\dots \}$,
$\N_0=\N \cup \{0\},$ and $p$ always 
denotes a prime number;
\item
The expression: $ f(\lambda,{\bf y},{\bf x}){\ll}_{{}_{{\bf y}}} g({\bf x})$
uniformly in ${ \bf x}\in X $ and ${\lambda}\in \Lambda$
means there exists $A=A({\bf y})>0$,
such that,   
$\forall {\bf x}\in X {\mbox { and }}\forall {\lambda}\in
{\Lambda}\quad |f(\lambda,{\bf y},{\bf x})|\leq Ag({\bf x}) $;
\item For any ${\bf x}=(x_1,..,x_n) \in \R^n$, we set   
$\|{\bf x}\|=\sqrt{x_1^2+..+x_n^2}$ and $|{\bf
  x}|=|x_1|+..+|x_n|$. We denote the canonical basis of $\R^n$ by 
$(\e_1,\dots,\e_n)$. The standard inert product on $\R^n$ is 
denoted by $\la .,. \ra$. We set also  
$\zerob=(0,\dots,0)$ and $\unb =(1,\dots,1)$;
\item 
We denote a vector in $\C^n$ $\s=(s_1,\dots,s_n)$, and write 
 $\s={\sigmab}+i{\taub},$
 where ${\sigmab } = (\sigma_1,\dots,\sigma_n)$ and  
 ${\taub }=(\tau_1,\dots,\tau_n)$ are the real resp. imaginary
 components of $\s$ (i.e. $\sigma_i=\Re(s_i)$ and 
$\tau_i=\Im(s_i)$ for all $i$). We also write 
$\la \x, \s\ra$ for $\sum_i x_i s_i$  if $\x \in \R^n, \s \in \C^n$;
\item Given $\alphab \in \N_0^n,$ we write $\X^{\alphab}$ for the
  monomial $X_1^{\alpha_1} \cdots X_n^{\alpha_n}$.  
For an analytic function  $h(\X)=\sum_{\alphab} a_{\alphab} \X^{\alphab}$, the set 
$supp(h):=\{\alphab \mid a_{\alphab} \neq 0\}$ is called the support of
$h$;
\item A function $f: \N^{n} \rightarrow \C$ is said to be 
  multiplicative if for all $m_1,\dots,m_n\in \N$ and 
  $m_1',\dots,m_n'\in \N$ satisfying 
$gcd\left(lcm\left(m_i\right),lcm\left(m_i'\right)\right)=1$ we have \\
$f\left(m_1 m_1',\dots , m_n m_n'\right)=f\left(m_1,\dots , m_n
\right) . f\left(m_1',\dots ,m_n'\right)$;
\item A polynomial $P\in \R [X_1,\dots, X_n]$ of degree $d$ is said
  to be elliptic if  its homogenuous part of highest  degree $P_d$ satisfies:
$\forall \x \in \R^n_+ \setminus \{\zerob\}$, $P_d(\x)>0$;
\item Let $F$ be a meromorphic function on a domain ${\cal D}$ of
  $\C^n$ and let ${\cal S}$ be the support of its polar divisor. 
$F$ is said to be of moderate growth if there exists 
$a,b>0$ such that $\forall \delta >0$,  
 $F(\s) \ll_{\sigmab,\delta} 1+|\tau|^{a|\sigmab|+b}$ 
uniformly in $s=\sigmab+i\taub\in {\cal D}$ 
verifying $d(\s, {\cal S})\geq \delta $. 
\ee
\subsection{Preliminaries from convex analysis}
For the reader's convenience, some classical notions from convex analysis 
that will be used throughout the article are assembled
here.\par
Let $A=\{\alphab^1,\dots,\alphab^q\}$ be a finite subset of 
  $\R^n$. 
\be 
\item The convex hull of $A$ is 
$conv(A):=\{\sum_{i=1}^q \lambda_i \alphab^i
\mid (\lambda_1,..,\lambda_q) \in \R_+^q {\mbox { and }} 
\sum_{i=1}^q \lambda_i =1\}$ and its interior is  
$conv^*(A):=\{\sum_{i=1}^q \lambda_i \alphab^i
\mid (\lambda_1,..,\lambda_q)\in \R_+^{*q} {\mbox { and }}
\sum_{i=1}^q \lambda_i =1\}$;
\item The convex cone of $A$ is $con (A):=\{\sum_{i=1}^q \lambda_i \alphab^i
\mid (\lambda_1,\dots,\lambda_q)\in \R_+^q \}$ and its  
(relative) interior is  
$con^* (A):=\{\sum_{i=1}^q \lambda_i \alphab^i
\mid (\lambda_1,\dots,\lambda_q)\in \R_+^{*q}\}$.
\ee 
Let $\Sigma $ be the set (or the boundary of the set) 
$\{ \x \in \R_+^n \mid \la \betab , \x \ra \geq 1 ~\forall \betab \in
I\}$ where $I$ is a finite (nonempty) subset of $\R_+^{n}\setminus
\{\zerob\}$. 
$\Sigma$ is a convex polyhedron of $\R_+^n\setminus
\{\zerob\}$.
\be 
\item Let $\a \in \R_+^n \setminus \{\zerob\}$, we define
  $m(\a):=\inf_{\x \in \Sigma } \la \a ,\x\ra$ 
and the face of $\Sigma $ with polar vector $\a$  (or the
first meet locus of $\a$)  as   
${\cal F}(\Sigma ) (\a)=\{\x \in \Sigma \mid \la \a , \x\ra =m(\a)\}$; 
\item The faces  of $\Sigma$ are the sets ${\cal F}(\Sigma ) (\a)$ 
$(\a \in \R_+^n \setminus \{\zerob\})$.
A facet  of $\Sigma $ is a   face  of maximal dimension;
\item
Let $F$ be a face of $\Sigma$.
The cone  $pol(F):=\{\a \in \R_+^n \setminus \{\zerob\} \mid
F={\cal F}(\Sigma ) (\a)\}$ is called the polar cone associated to $F$  
and its elements are called polar vectors of $F$.
\item We define the {\it index} of  $\Sigma$ 
by $\iota (\Sigma):=\min\{|\alphab |;~
\alphab \in \Sigma \}$. It is  clear that \\
${\cal F}(\Sigma ) (\unb)=\{\x \in \Sigma;~ |\x| = \iota(\Sigma)\}$.
\ee 
\subsection{Construction of the volume constant}
{\bf The Sargos  constant (\cite{sargosthese}, chap 3, \S 1.3) :}\\
Let $P(\X)= \sum_{\alphab \in supp(P)} a_{\alphab} \X^{\alphab}$ be 
a generalized polynomial; i.e. $supp(P)$ is a finite subset of 
$\R_+^n$ (and not necessarly of $\N_0^n$). We suppose that $P$ has 
positive coefficients and that it depends on all the variables $X_1,\dots,X_n$.
We denotes by 
${\cal E}^{\infty}(P):= \left(conv(supp(P))-\R_+^n\right)$ its
Newton polyhedron at infinity. 
Let $G_0$ be the smallest face of ${\cal E}^{\infty}(P)$ which meets
the diagonal $\Delta =\R_+ \unb $.
We denote by $\sigma_0=\sigma_0(P)$ the unique positive real number
$t$   that satisfies  $t^{-1}\unb \in G_0.$ We also set $\rho_0 =\rho_0 (P):=codim G_0$.\\ 
By a permutation of coordinates one can suppose that $\oplus_{i=1}^{\rho_0}\R \e_i
\oplus \overrightarrow{G_0}=\R^n$ and that \\
$\{\e_i \mid G_0=G_0-\R_+ \e_i\}=\{\e_{m+1},\dots,\e_n\}$.\\
Let $\lambdab_1,\dots,\lambdab_N$ be the polar vectors of the facets
of ${\cal E}^{\infty}(P)$ which meet $\Delta$.
Set $P_{G_0}(X)=\sum_{\alphab \in G_0} a_{\alphab} \X^{\alphab}$ and 
$\Lambda =Conv\{\zerob,\lambdab_1,\dots,\lambdab_N, \e_{\rho_0+1},\dots,
\e_n\}$.
\begin{Definition} 
The Sargos constant associated to $P$ is:\\
$
A_0(P):= n!~Vol(\Lambda)~\int_{[1,+\infty[^{n-m}}\left(\int_{\R_+^{n-\rho_0}}P_{G_0}^{-\sigma_0}
(\unb,\x, \y) ~d\x\right) d\y >0.
$
\end{Definition}
In (\cite{sargosthese}, chap 3), P. Sargos proved 
the following important result:\\
{\bf Theorem (\cite{sargosthese}, chap 3, th. 1.6):} \label{sargosth}
{\it 
Let $P$ be a generalized polynomial as above. We set 
$Y(P;s):=\int_{[1,+\infty[^n } P(\x)^{-s} ~d\x$. 
The abscissa of convergence of $Y(P;s)$ is $\sigma_0=\sigma_0(P)$. 
Moreover $s\mapsto Y(P;s)$ has a meromorphic continuation to $\C$,    
$\sigma_0$ is indeed a pole of $Y(P;s)$ of order $\rho_0$ and   
$\ds Y(P;s)\sim_{s\rightarrow \sigma_0}
 A_0(P) ~(s-\sigma_0)^{-\rho_0}$.
}\par

\medskip
If $P$ is elliptic the previous result can be sharpened as follows:
\begin{Proposition} [\cite{mahler}] \label{sargoselliptique}
Let $P \in \R [X_1,\dots,X_n]$ be an elliptic polynomial of  
degree $d\geq 1$. Denote by $P_d$ its homogeneous part of greater
degree. Then  
$\sigma_0 =\frac{n}{d}$, $\rho_0 =1$ and the Sargos constant
associated to $P$ is:
$\ds A_0(P)=\frac{1}{d} \int_{\Sp^{n-1} \cap \R_+^{n}} P_d^{-n/d}(\v)
d\sigma(\v)$ where $\Sp^{n-1}$ is the unit sphere of $\R^n$ and
$d\sigma $ is induced Lebesgue measure.
\end{Proposition}
{\bf Construction of the volume constant:}\\
Let $I$ be a finite subset of $\R_+^r \setminus \{\zerob\}$, 
$\u=(u(\betab)_{\betab \in I})$ a finite sequence of elements of $\N,$
and $\b =(b_1,\dots,b_r)\in \R_+^{*r}$.
To this data, we associate the generalized polynomial  
$P_{(I;\u;\b)}$ with $q:=\sum_{\betab \in I} u(\betab)$ variables, 
in the following way:\\
We define $\alphab^1,\dots,\alphab^q$ by:
$\{\alphab^i \mid i=1,\dots,q\}=I$ and
$\forall \betab \in I$ 
$\#\{i\in \{1,..,q\} \mid \alphab^i=\betab \}=u(\betab)$ (i.e. the
family $(\alphab^i)$ is obtained by repeating  each   $\betab$   $u(\betab)$ times).\\
We define the vectors $\gammab^1,\dots,\gammab^r$ 
of $\R_+^q$ by: $\forall i=1,\dots,q$ and  $\forall k=1,\dots,r$, $\alpha_k^i=\gamma_i^k$.\\
We set finally $P_{(I;\u;\b)}(\X):=\sum_{i=1}^r b_i
\X^{\gammab^i}$.\\
We define the {\it volume constant  
  associated to}  $I$, $\u$, $\b$ by
$\ds A_0(I;\u;\b):=A_0\left(P_{(I;\u;\b)}\right)>0$.

\subsection{Some important constants: mixed exponents and volume}
Let  
$P(\X)=b_1 \X^{\gammab^1}+\dots+b_r
\X^{\gammab^r}\in \R_+[X_1,\dots,X_n]$. We set $\b =(b_1,\dots,b_r)\in
\R_+^{*r}$.\\ 
Let $\T=(I,\u)$ where $I$ is a finite subset of 
$\R_+^n\setminus \{0\}$ and $\u =\left(u(\betab)\right)_{\betab \in I}$
a vector of positive integers.\par
We associate to $\T$ and $P$ the following (mixed) objects which will 
play an important role in the sequel of this paper:
\be
\item $n$ elements $\alphab^1,\dots,\alphab^n$ of
  $\N_0^r$ defined by:
$\alpha^i_j =\gamma^j_i$ $\forall i=1,\dots,n$ and  $\forall
j=1,\dots,r$;
\item $\mu(\T;P;\betab):=\sum_{i=1}^n \beta_i \alphab^i$ for all  
$\betab \in I$ \ \ and \ \ $I_{\T,P}=\{\mu(\T;P;\betab) \mid \betab \in
I \}$; 
\item $\u_{\T,P}=\left(u_{\T,P}(\etab)\right)_{\etab \in I_{\T,P}}$ \ \ where \ \  
$u_{\T,P}(\etab)=\sum_{\{\betab \in I;
  ~\mu(\T;P;\betab)=\etab\}} u(\betab)$  
$\forall \etab \in I_{\T;P}$;
\item 
$\rho_0 (\T):= \sum_{\betab \in I} u(\betab) -rank\left(I \right)
+1$ and  $\rho_0 (\T;P):= \sum_{\etab \in I_{\T,P}} u_{\T,P}(\etab)
  -rank(I_{\T,P}) +1$.
\ee
We will also use the  following easy to check {\bf remark:} \\
{\it $\rho_0 (\T;P)=\rho_0(\T)$
if $rank \,\left(supp(P)\right)=n$. In particular this is true if $P$ is elliptic.}\par

We define finally the {\it mixed volume} constant by:
$$
A_0(\T;P)=A_0(I_{\T,P};\u_{\T,P};\b)
$$
where $A_0(I_{\T,P};\u_{f,P};\b)>0$ is the volume constant 
 associated to $I_{\T,P};\u_{\T,P}$ and $\b$.

\section{Statements of Main Results}
\subsection{Functions of finite type: Definition and  
examples}
Let $f:\N^{n}\mapsto \C$ be a  function. We define formally 
${\cal M}(f;\s):=\sum_{\m \in \N^{n}} \frac{f(m_1,\dots,m_n)}{m^{s_1}\dots
  m_n^{s_n}}$.\\ 
Typically, though not always, the series that are of interest  in number theory are created    by expanding out an  
Euler   product in one or more variables, in which case  the function
$\m \to  f(\m)$  will be multiplicative. 
Several works (see for example \cite{kurokawa},
\cite{morozeuler}, \cite{bretechecompo}, \cite{eulerprod}...) then indicate that  
${\cal M}(f;\s)$ should satisfy the following general property:
\begin{Definition}\label{fonctiontf}
A function $f:\N^{n}\rightarrow \C$ is said to be of finite type if there exists \\
a finite subset $I$ of $\R_+^n\setminus \{\zerob\}$ such that if we
set for all $\delta \in \R$: 
$$V(\delta):=\{\s \in \C^n \mid \la \betab, \Re \s \ra> \delta ~\forall
\betab \in I\}\cap \{\s \in \C^n \mid \Re(s_i) >0 \forall
i=1,\dots,n\},$$  then :\\
(i) ${\cal M}(f;\s)$ converges absolutely in $V(1)$. \\
(ii)  there exist $\delta <1$ such that ${\cal M}(f;\s)$ can
be continued to the set $V(\delta)$ as a 
meromorphic function  with moderate growth in $\Im s$ and polar
divisor 
$\bigcup_{\betab \in I} \{\la \betab, \s \ra = 1\}$.
\end{Definition}   
We define $\Sigma_f = \R_+^{n} \cap \partial \left(\{\x \in \R^n \mid \la
\betab, \x\ra > 1 | \forall \betab \in I\}\right) \subset  \partial V(1)$.\\
For each $\cb \in \Sigma_f \cap \R_+^{*n}$, we associate the pair $\T_{\cb}
:=(I_{\cb},\u ),$ where 
$I_{\cb}:=\{\betab \in I \mid \la \betab , \cb \ra=1\}$ and $u =
\big(u(\betab)\big)_{\betab \in I_{\cb}}$ is a vector of positive
integers defined by the following property: 
there exists   
 $\eps_0 >0$  such that 
$\ds \s\mapsto H_{\cb}(f; \s):=\big(\prod_{\betab \in I_{\cb}} {\left(\la
    \s,\betab\ra \right)}^{u(\betab)}\big)~{\cal M}(f;\cb+\s)$ 
has a {\it holomorphic} continuation  to the set 
$\{ \s \in \C^n \mid \sigma_i >-\eps_0 ~\forall i\},$ and does not
  vanish identically along $\la \s,\betab\ra = 0$. 
Thus, $u(\beta)$ equals  the order of the pole of
  ${\cal M}(f;\cb+\s)$ along $\{\la \s,\betab\ra
=0\}$ at $\s = \zerob.$\par
We call the pair $\T_{\cb}$  the {\it polar type} of $f$ at $\cb$.\par 
{\bf Examples of functions of finite type:} 
\be
\item If $f$ is a monomial possibly twisted by additive characters
  then $f$ is clearly of finite type. We call this case the classical case;
\item Uniform multiplicative functions (see remark 5 \S 3.2), among
  which are characteristic functions of multiplicative sets;
\item  The functions $f$ where the Dirichlet series ${\cal
    M}(f, s)$ belongs to the class of Dirichlet series with Euler
  product  studied by  N. Kurokawa \cite{kurokawa}
and B.Z. Moroz \cite{morozeuler}.
\ee
\subsection{Results about mixed zeta functions}
Let $P(X)\in \R_+[X_1,\dots,X_n]$ be a homogeneous polynomial with
positive coefficients which depend on all the variables 
$X_1,\dots,X_n$. We denote by $d\geq 1$ its degree.
Let $f$ be a function of finite type (see definition
\ref{fonctiontf}.)
Set
$\ds Z(f;P;s):=\sum_{\m \in \N^{n}}
\frac{f(m_1,\dots,m_n)}{P(m_1,\dots,m_n)^{s/d}}$.\\ 
Recall (see \S 2.2) that   ${\cal F}(\Sigma_f)(\unb)$ denotes the face of
$\Sigma_f$ whose polar vector is the diagonal vector $\unb$. For simplicity, the index $\iota(\Sigma_f)$ is denoted 
$\iota(f)$ in the following.\\
We now define the mixed data as follows:\\  
$$\Sigma_f (P) := \Sigma_f   \cap 
con^*\left(supp(P)\right), \qquad 
\Sigma_f (P;\unb):= {\cal F}(\Sigma_f)(\unb) \cap con^*\left(supp(P)\right),$$
  and the two indices $$\iota (f;P)=\inf \{|\cb| : 
\cb
\in \Sigma_f (P)\}, 
\qquad \rho (f;P)= min\{\rho_0(\T_{\cb};P) : \cb \in \Sigma_f (P){\mbox
  { and }} |\cb|=\iota(f;P)\}. $$ 
The main results of this paper are the following two theorems.
\begin{Theoreme}\label{analytic}
Assume that $\Sigma_f (P) \neq \emptyset$.
Then $s \mapsto Z(f;P;s)$ is {\it holomorphic} in the half-plane 
$\{s: \sigma > \iota (f;P)\}$ and there exists $ \eta >0$ such that   
 $s \mapsto Z(f;P;s)$ has a meromorphic continuation with moderate
 growth to the half-plane $\{\sigma > \iota (f;P) - \eta\}$ with at most
one pole at $s=\iota (f;P).$   If $\iota(f;P)$ is a pole, then its order is at most $\rho (f;P)$.    
\end{Theoreme}
{\bf Remark.} In general, one should not expect $\Sigma_f(P;\bold 1)$ to be nonempty when $\Sigma_f(P) \neq \emptyset.$
As noted  in \S 2.2, the function $\cb \in \Sigma_f \to |\cb |$ assumes its minimal value
$\iota (f)$ on the face ${\cal F}(\Sigma_f)(\unb)$ of $\Sigma_f$. Thus, when $\Sigma_f(P;\bold 1) =
\emptyset,$ it follows that $\iota(f) < \iota(f;P).$ However, if $\Sigma_f(P;\bold 1) \neq \emptyset,$ then $\iota(f) =
\iota(f;P).$ In this case,  Theorem
\ref{analytic}   says that 
$Z(f;P;s)$ must always remain  analytic  up to, at least, the index of
$\Sigma_f.$ \qed \par 
Although Theorem 1 is quite general, it is also not as explicit as one would like in practice. The following theorem 
gives, under some natural assumptions, a more explicit description of the domain of analyticity, as well as  the
principal part of $Z(f;P;s)$ at its first pole:
\begin{Theoreme}\label{partieprincipale} 
Assume that there exists $\cb \in \Sigma_f (P;\unb)$ such
that:
\be
\item $\unb \in con^*(I_{\cb})$; 
\item there exists a function $L$ such that: 
$\ds {\cal M}(f;\s)= L\left((\la \betab,\s\ra)_{\betab \in
 I_{\cb}}\right)$.
\item $H_{\cb}(f;\bold 0) \neq 0.$
\ee
Then  $s=\iota(f)$ is indeed a pole of order $\rho (f;P)$ and  
$$Z(f;P;s)\sim_{s\rightarrow \iota(f)}
\frac{C(f;P)}{\left(s-\iota(f)\right)^{\rho (f;P)}}\,, \qquad {\mbox { where }}
C(f;P):= H_{\cb}(f;\zerob) d^{\rho (f;P)} A_0(\T_{\cb}, P)\neq 0.$$ 
\end{Theoreme}
Let me now give some additional remarks:\\
{\bf Remark 1.} If $P$ is elliptic, then
  $con^*\left(supp(P)\right)=\R_+^{*n}$. Consequently  
$\Sigma_f (P;\unb) = {\cal F}(\Sigma_f)(\unb)$.\par
{\bf Remark 2} From the definition of 
$I_{\cb}$ it follows that the smallest face which contains $\cb$ is 
$\bigcap_{\betab \in I_{\cb}} {\cal F}(\Sigma_f)(\betab).$ A well known
result of convex analysis then implies  that the set of its polar vectors 
is $con^*(I_{\cb})$. So assumption (1) of theorem
\ref{partieprincipale} is equivalent to:\\ 
$(1)'$ $\quad {\cal F}(\Sigma_f)(\unb)$ is the smallest face of
$\Sigma_f$ which contains $\cb$.\par
{\bf Remark 3.} The assumption (2) of theorem \ref{partieprincipale} is always 
satisfied if $rank \left(I_{\cb}\right)=n$, in particular 
if $\cb$ is a vertex of $\Sigma_f$.\par
{\bf Remark 4.}
In the classical case where $f(\m)=\m^{\mub}$ is a monomial, we
  have ${\cal M}(f;\s)=\prod_{i=1}^n \zeta(s_i -\mu_i)$. Thus it is
  clear that $f$ is of finite type. Moreover, if we set 
$I=\{\frac{1}{1+\mu_i} \e_i \mid i=1,\dots,n\}$, then the domain of
convergence of ${\cal M}(f;\s)$ equals ${\cal D}_f=\{\s \in \C^n; ~\la
\betab , \s\ra >1~\forall \betab \in I\}.$  
In particular  it is clear that  
$\cb=(1+\mu_1,\dots,1+\mu_n) \in     {\cal
  F}(\Sigma_f)(\unb)$  and that the polar type of $f$ in $\cb$ is 
$\T_{\cb}=(I;\u)$ where $u(\betab) =1$ for all $\betab \in I$. 
In addition,   conditions  (1), (2) and (3) of Theorem
\ref{partieprincipale} are easily verified.  
Thus, if the polynomial $P$ satisfies $\cb \in
 con^*\left(supp(P)\right),$  then the conclusion of Theorem 2
 follows.  \\We will  check this if  $\mub = \zerob$ 
and $P$ is elliptic. In this case, it is clear that   
  $\iota(f)=|\cb|=n$    and  
$\rho(f;P)=\rho_0(\T_{\cb};P)=\sum_{\betab \in I} u(\betab)
-rank(I) +1=1$. Thus,   
$$Z(P;s)=\sum_{\m \in \N^{n}} P^{-s/d}(\m)=Z(\unb;P; s)
\sim_{s\rightarrow \iota(f)}
\frac{A_0(\T_{\cb},P)}{\left(s-n \right)}.$$
Moreover, it is easy to see that in this case $A_0(\T_{\cb},P)=A_0(P)$ the Sargos
constant (see \S 2.3). So we find  Mahler's result \cite{mahler}
without the factor $1/d$ (due to the  normalization $s \to
s/d$). \par
{\bf Remark 5.}
Let $f:\N^n\rightarrow \C$ be a multiplicative function. We suppose that
$f$ is uniform i.e. that there exist two constants $C, M>0$ 
such that for all prime $p$ and all $\nub \in \N_0^n$ 
$f(p^{\nu_1},\dots,p^{\nu_1})$ is independent of $p$ and verifies
$f(p^{\nu_1},\dots,p^{\nu_n}) \leq C (1+|\nub|^M)$.\\
Set $S^*(f)=\{\nub \in \N_0^n\setminus
\{\zerob\} \mid f(p^{\nu_1},\dots, p^{\nu_n}) \neq 0\}$. Denote by 
$E(f):={\cal E}(S^*(f))= conv \left( S^*(f) +\R_+^n\right)$ 
its Newton polyhedron. We assume that ${\cal
  F}\left(E(f)\right),$ its smallest face which meets the diagonal, is compact.
Let $\cb \in \R_+^{*n}$ be a normalized polar
  vector of ${\cal F}\left(E(f)\right)$ 
(i.e. ${\cal F}\left(E(f)\right)=E(f)\cap \{\x \in \R^n \mid \la \cb ,
\x\ra =1\}$). Then it follows from \cite{essouametrics}
that $f$ is a function of finite type, that
the polar type of $f$ at $\cb$ is $\T=(I_{\cb};\u),$ 
where $I_{\cb}:={\cal F}\left(E(f)\right)\cap S^*(f),$
and 
$\u= \left(f(p^{\beta_1},\dots, p^{\beta_n})\right)_{\betab \in I}$. 
Moreover if $P$ is elliptic then $\cb \in
\Sigma_f(P;\unb)$ and:
\be 
\item Assumption (1) of theorem \ref{partieprincipale} 
is satisfied by definition of $\cb$; 
\item Assumption (2) of theorem \ref{partieprincipale} 
is satisfied if, for example, 
$rank({\cal F}\left(E(f)\right))=rank(S^*(f))$.
\item Assumption (3) of
  theorem is  \ref{partieprincipale}
also satisfied by  (\cite{eulerprod}, Theorem 6). 
\ee
In particular if  ${\cal F}\left(E(f)\right)$ is a 
compact facet of the polyhedron $E(f)$, then all the assumptions of theorem
\ref{partieprincipale} are satisfied.  
\subsection{General counting functions and 
  lattice points problems}
By a simple adaptation of a standard tauberian argument of Landau (see
for example \cite{essouabrismf}, Prop. 3.1)), 
we obtain, with the notations of theorems \ref{analytic} and
\ref{partieprincipale}, the following:
\begin{Corollaire}\label{applicationneuve}
  Let $f:\N^n\rightarrow \R_+$ be a function of finite type and $P\in
  \R_+[X_1,\dots,X_n]$ be a homogeneous polynomial of degree $d\geq 1$. Define: 
  $$N(f;P;t):=\sum_{\{\m \in \N^{n}; ~P(\m)^{1/d}\leq t\}} f(m_1,\dots,m_n), \qquad t > 0. $$ 
\indent Assume that 
$\Sigma_f (P) \neq \emptyset$. Then $\ds N(f;P;t) = O\left(t^{\iota(f;P)}~(\log t)^{\rho
  (f;P)-1}\right) {\mbox { as }} t\rightarrow \infty$.\par
\quad Now assume that there exists $\cb \in \Sigma_f  (P;\unb)$
satisfying the conditions of theorem \ref{partieprincipale}. 
Then there exist $\theta >0$ and 
a polynomial $Q\neq 0$ of degree $\rho (f;P) -1$ such that as $t\rightarrow
\infty$ 
$$
N(f;P;t)= t^{\iota(f)} Q(\log t)+O\left(t^{\iota(f)-\theta}\right)
= C_0(f;P) ~ t^{\iota(f)} (\log t)^{\rho (f,P) -1} ~
\left(1+O\left((\log t)^{-1}\right)\right),
$$
where   
$$C_0(f;P):=\frac{C(f;P)}{\iota(f) ~ (\rho (f;P)-1)!}
=\frac{H_{\cb}(f;\zerob) d^{\rho (f;P)}
 A_0(\T_{\cb},P)}{\iota(f) ~ (\rho (f;P)-1)!} >0. $$
 \end{Corollaire}
{\bf Remark:} Corollary  \ref{applicationneuve} estimates counting 
functions associated to a large class 
of multivariable  arithmetic functions and polynomials. 
To our knowledge, the only comparable result 
 is due to de la Bretèche \cite{bretechecompo} who proved estimates for counting 
functions of the form $N_\infty(f;t):=\sum_{\{\m \in \N^{n}; ~\max_i m_i\leq t\}}
f(m_1,\dots,m_n)$ for functions $f$ satisfying conditions similar to
those in Definition \ref{fonctiontf}.  So corollary
\ref{applicationneuve} can be viewed as an extension of his result to  a class of generalized heights determined by a
homogeneous $P \in  \R_+[X_1,\dots,X_n].$\cqfd \par
 
We apply this Corollary to some cases of arithmetic interest. The first  example is motivated by the work of K.
Matsumoto  and Y. Tanigawa   \cite{matsumotobordeaux}. 
Consider $n$ arithmetical functions $f_j: \N \rightarrow \C$
$(j=1,\dots,n)$. 
We assume that for each $j$ there exist $c_j >0$, $u_j \in \N$
and $\eta_j >0$ such that the zeta function  
$Z(f_j;s):=\sum_{m=1}^{\infty} \frac{f_j(m)}{m^s}$.
converges absolutely in the half-plane $\{\sigma
>c_j\}$ and has a
meromorphic continuation with moderate growth to $\{\sigma >c_j
-\eta\}$  whose only  pole occurs  at $s=c_j$ with order $u_j$.\par
Let $P(X)\in \R[X_1,\dots,X_n]$ be elliptic of degree
$d\geq 1$.\par
Then, for   $\Re s$ sufficiently large, the series:
$\ds Z(\f;P;s):= \sum_{\m \in \N^n} \frac{f_1(m_1)\dots
  f_n(m_n)}{P^{s/d}(m_1,\dots,m_n)}\,,$\\
converges absolutely. For $t>0$ define: 
$\ds N(\f;P;t):=\sum_{\{\m \in \N^n; ~P^{1/d}(\m)\leq t\}} 
f_1(m_1)\dots f_n(m_n)$.\par
Set $\cb=(c_1,\dots,c_n)$, $c:=|\cb|= c_1+\dots+c_n$ and 
$\rho_0:=\sum_{j=1}^{n} u_j -n+1$.\\ 
Also, set $I:=\{c_j^{-1} \e_j \mid j=1,\dots,n\}$,
$\u=\left(u(\betab)\right)_{\betab \in I}\,,$ where $u(c_j^{-1}
\e_j)=u_j$ for all $j$, $\T_0 =(I;\u),$ and $A_0(\T_0;P)$ the mixed
volume constant associated to $\T_0$ and $P$.\\
Finally, for each $j=1,\dots,n,$ set $A_j=\lim_{s\rightarrow 0} s^{u_j}
Z(f_j;P;c_j+s) \neq 0$ and $A:=\prod_{j=1}^n A_j$.  
 
\begin{Corollaire}\label{mats+} 
There exists $\eta >0$ such that 
$s\mapsto Z(\f;P;s)$ has a meromorphic continuation with moderate growth
to $\{\sigma > c -\eta\}$ whose only pole is $s=c$ with order equal to 
$\rho_0$. Moreover:
$\ds Z(\f;P;s)\sim_{s\rightarrow c} A~ d^{\rho_0 }
A_0(\T_0,P) ~\left(s-c\right)^{-\rho_0 }$.\par
If we assume in addition that $ f_1(m_1)\dots f_n(m_n) \geq 0$ then 
there exists $\delta >0$  and 
a polynomial $Q\neq 0$ of degree $\rho_0 -1$ such that
$$
N(\f;P;t)= t^{c} Q(\log t)+O(t^{c-\delta})
= C~t^c \log^{\rho_0-1} t \cdot \left(1 + O(\frac 1{\log t})\right)
{\mbox { as }} t\rightarrow \infty,$$ where  
$\ds C= \frac{A~ d^{\rho_0 }
A_0(\T_0,P)}{c~(\rho_0 -1)!}$.
\end{Corollaire}
{\bf Remark.}   In fact,  \cite{matsumotobordeaux} only studied  the 
    meromorphic continuation of the multiple zeta function  
$\sum_{\m \in \N^n} \frac{f_1(m_1)\dots f_n(m_n)}{m_1^{s_1} \dots
    (m_1+\dots+m_n)^{s_n}}$. The point of Corollary \ref{mats+} is that it 
extends [ibid.] by noting that   Corollary \ref{applicationneuve} can
    be applied to deduce precise information about the principal part
    of any mixed zeta function
$Z(\f;P;s)$ at its first
    pole,  as well as an explicit asymptotic of $N(\f;P;t)$  
whenever   $\f(\m) := f_1(m_1)\dots f_n(m_n) \geq 0$. \par
A particular application of this corollary is therefore the
    following. Let $k\geq 2$ be a positive integer and let 
$B_k$ denote the set of lattice points
$\m \in \N^{n}$ such that each $m_i$ is $``k-$free" (i.e. $v_p(m_i)\leq
    k-1$ for any prime $p$) and set ${\unb}_{B_k}$ to be
the characteristic function of $B_k$.  
\begin{Corollaire}\label{squarefree}
Let $P\in \R_+[X_1,\dots,X_n]$ be elliptic of degree 
$d\geq 1$. Set 
$$N ({\unb}_{B_k}; P;t):=\#\{ \m \in \N^{n} \mid P^{1/d}(m_1,\dots,m_n)\leq t
{\mbox { and }} 
 m_i {\mbox { is }}k-{\mbox {free }} \forall i~\}.$$  Then there exists $\theta
>0$ such that:
$N ({\unb}_{B_k};P;t)=C_n(B_k, P)~ t^{n} +O\left(t^{n -\theta}\right),$
where 
$$
C_n(B_k, P):=  \left(\frac{1}{\zeta(k)}\right)^n \,\cdot \frac{1}{n} \,
\int_{\Sp^{n-1} \cap \R_+^n} P_d^{-n/d}(\v)~ d\sigma(\v).$$
\end{Corollaire}
{\bf Remark.} $B_2$ is the set of lattice points
$\m \in \N^{n}$ such that each $m_i$ is squarefree. In this case 
Corollary \ref{squarefree} gives  
$ C_n(B_2, P):= \left(\frac{6}{\pi^2}\right)^n \frac{1}{n}
\int_{\Sp^{n-1} \cap \R_+^n} P_d^{-n/d}(\v)~ d\sigma(\v)$.\par
A second application of Corollary \ref{applicationneuve} concerns the subset $D_k$ of $\m \in \N^n$ such
that $m_1\cdots m_n$ is $k-$free. This set is not 
  the $n-$fold product of a multiplicative subset of $\N$ whose asymptotic is 
standard. So, we should not expect the constant in the main term to have as simple a structure as in
Corollary \ref{squarefree}. 
\begin{Corollaire}\label{squarefreeuniform}
Let $P\in \R_+[X_1,\dots,X_n]$ be elliptic of degree 
$d\geq 1$. Let $k\geq 2$ be a positive integer and set  
$$N ({\unb}_{D_k};P;t) = \# \{\m = (m_1,\dots,m_n) \in \N^n :  P^{1/d}(\m) \le t  \ \text{and  $m_1\dots m_n$  is
$k-$free}\}.$$  
 Then there exists $\theta >0$ such that:
$N ({\unb}_{D_k};P;t) = C_n\left(D_k, P\right)~ t^{n} +O\left(t^{n -\theta}\right),$
where  
$$
C_n\left(D_k, P\right):= \left[\prod_p\left(1-\frac{1}{p}\right)^n
\left(1+\sum_{l=1}^{k-1}\frac{{n+l-1\choose l}}{p^l}\right)\right]~\frac{1}{n}
\int_{\Sp^{n-1} \cap \R_+^n} P_d^{-n/d}(\v)~ d\sigma(\v).$$ 
\end{Corollaire}
 
A fourth class of examples is restricted to $n-$tuples of the 
primes ${\cal P}.$ Let $H \in \Z_+[X_1,\dots,X_n]$ be a
polynomial.  For any positive integer $k$, set $V_k (H):=\{ (p_1,\dots,p_n)\in
{\cal P}^n \mid k=H(p_1,\dots,p_n)\}$.  The following corollary   gives   some
results about $V_k (H)$ on average, assuming (RH). 
Set for all $t>0$:
\be
\item $N^{(0)}({\cal P}^n;H;t):= \sum_{k\leq t} V_k (H)= 
\sum_{\{(p_1,\dots,p_n)\in {\cal P}^n; ~H(p_1,\dots,p_n)\leq t\}} 1$,
\item 
$N^{(1)}({\cal P}^n;H;t):=  
\sum_{\{(p_1,\dots,p_n)\in {\cal P}^n; ~H(p_1,\dots,p_n)\leq t\}} \log
p_1 \dots \log p_n$  
\item 
$N^{(2)}({\cal P}^n;H;t):= 
\sum_{\{(m_1,\dots,m_n)\in \N^n; ~H(m_1,\dots,m_n)\leq t\}}
\Lambda (m_1) \dots \Lambda (m_n)$ where $\Lambda$ denotes    
Mangoldt's function.
\ee
In order to obtain some results toward Goldbach's conjecture on average, the
functions $N^{(1)}({\cal P}^n;H;t)$ and $N^{(2)}({\cal P}^n;H;t)$ have
been studied by several papers (\cite{fujiiacta}, \cite{fujiipja}, ..) in the
particular case $H(X_1,X_2)=X_1+X_2$. In such work, the important point has been to 
get as good an error term as possible. One can, however, without much additional effort, 
derive from corollary \ref{mats+} above an explicit main term 
for a large class of nonlinear polynomials $H$ as follows.
\begin{Corollaire}\label{goldba}
Let $H \in \R_+[X_1,\dots,X_n]$ be an elliptic polynomial of degree
$d\geq 1$.  
Assuming $(RH)$, there exists $\delta >0$ such that 
for each $i\in \{1,2\}$:\\
$N^{(i)}({\cal P}^n;H;t)=
C~t^{n/d} +O( t^{n/d -\delta})\,,$ where 
$C=\frac{1}{n} 
\int_{\Sp^{n-1} \cap \R_+^n} H^{-n/d}(\v)~ d\sigma(\v)$.
\end{Corollaire}
{\bf Remark.} By using corollary \ref{goldba}, it is easy to see that 
for any elliptic polynomial $H$ of degree $d\geq 1$, there exist
$\alpha ,\beta >0$ such that for all $t>0$:
$$\alpha ~ t^{n/d}(\log t)^{-n} \leq N^{(0)}({\cal P}^n;H;t):= 
\sum_{k\leq t} V_k (H) \leq \beta ~t^{n/d}.$$
\section{Proofs  of Theorems 1 and 2 and their corollaries}
\subsection{Fundamental lemmas and their proofs}
In this section we will gather some important lemmas for the proofs 
of our theorems. Lemma \ref{majoration}, rather elementary, will enable us to justify 
convergences of the integrals in our proofs. Lemmas \ref{fonda1} and
\ref{fonda3} are the heart of our method and constitute 
the most important technical part of this work. 

\begin{Lemme}\label{majoration}
Let $\mub^1,\dots,\mub^r$ be vectors of $\R^n$ and let  $l\in \R$.
Set for all  $\tau \in \R$:
$$F_n(\tau):=\int_{\R^n} \prod_{i=1}^n (1+|y_i|)^l
e^{-\frac{\pi}{2}\left(\sum_{i=1}^n |y_i|+
\sum_{i=1}^n |\la \mub^i,\y \ra| + |\tau -\sum_{i=1}^n y_i 
- \sum_{i=1}^n \la \mub^i,\y \ra|\right)}~dy_1\dots dy_n$$ 
There exist $A=A(l,n), C = C(l,n) >0$ such that $\forall \tau \in \R$:
$ F_n(\tau) \leq C~(1+|\tau|)^A ~e^{-\frac{\pi}{2}|\tau|}$.
\end{Lemme}
{\bf Proof of the Lemma \ref{majoration}:}\\
Set  
$\psi_n(l;\tau):=\int_{\R^n} \prod_{i=1}^n (1+|y_i|)^l
e^{-\frac{\pi}{2}\left(\sum_{i=1}^n |y_i| + |\tau -\sum_{i=1}^n y_i| 
\right)}~dy_1\dots dy_n.$
Since  
$$ 
\sum_{i=1}^n |\la \mub^i,\y \ra| + \left|\tau -|\y|
- \sum_{i=1}^n \la \mub^i,\y \ra\right| 
\geq 
\left|\sum_{i=1}^n \la \mub^i,\y \ra + \tau -|\y|
- \sum_{i=1}^n \la \mub^i,\y \ra \right|
\geq 
\left|\tau - |\y|\right|.$$
It follows that $F_n(\tau) \leq \psi_n(l;\tau)$. 
So to prove the lemma it is sufficient to prove the same inequality for 
$\psi_n(l;\tau)$. We now proceed by induction on $n$.\par
$\bullet$ {\bf Case n=1:}\\
It is sufficient to prove the inequality for 
$\psi_1^{+}(l;\tau) =  \int_0^{+\infty} (1+y)^l 
e^{-\frac{\pi}{2}\left(y + |\tau -y | 
\right)}~dy$ and $\tau >1$. The inequalities for the other 
cases rise in a similar way.\\  
For $\tau >1$ we have:
\begin{eqnarray*}
\psi_1^{+}(l;\tau) &= & \int_0^{+\infty} (1+y)^l 
e^{-\frac{\pi}{2}\left(y + |\tau -y | 
\right)}~dy\\
& = & \int_0^{\tau} (1+y)^l 
e^{-\frac{\pi}{2}\left(y + |\tau -y | 
\right)}~dy + \int_{\tau}^{+\infty} (1+y)^l 
e^{-\frac{\pi}{2}\left(y + |\tau -y | 
\right)}~dy\\
& \ll_l & 
\tau^{l+1}+ e^{\frac{\pi}{2}\tau} \int_{\tau}^{+\infty} (1+y)^l 
e^{-\pi y } ~dy 
\ll_l 
\tau^{l+1}+ e^{-\frac{\pi}{2}\tau} \int_{1}^{+\infty} (t+\tau )^l 
e^{-\pi t } ~dt \\
& \ll_l &
\tau^{l+1}+ (\tau +1)^l e^{-\frac{\pi}{2}\tau} \int_{1}^{+\infty} (1+t)^l
e^{-\pi t } ~dt \ll_l (\tau +1)^{l+1} e^{-\frac{\pi}{2}\tau}.
\end{eqnarray*}
This proves the lemma for $n=1$.\\
$\bullet$ {\bf End of the induction:}\\
Let $n\geq 2$. Assume that the lemma is true for $n-1$.
Thus there exists $A>0$ such that\\ 
$\psi_{n-1}(l;\tau) \ll_l (1+|\tau|)^A ~e^{-\frac{\pi}{2}|\tau|}$ 
($\tau \in \R$). 
It follows that we have uniformly in $\tau \in \R$:
\begin{eqnarray*}
\psi_n(l;\tau)&:=& \int_{\R^n} \prod_{i=1}^n (1+|y_i|)^l
e^{-\frac{\pi}{2}\left(\sum_{i=1}^n |y_i| + |\tau -\sum_{i=1}^n y_i| 
\right)}~dy_1\dots dy_n \\
& = & \int_\R \psi_{n-1}(l;\tau -y_n)~(1+|y_n|)^l
e^{-\frac{\pi}{2}|y_n|}~dy_n\\
& \ll_{l,n} &  \int_\R ~(1+|\tau -y_n|)^A e^{-\frac{\pi}{2}|\tau - y_n|}(1+|y_n|)^l
e^{-\frac{\pi}{2}|y_n|}~dy_n \quad ({\mbox { by the induction hypothesis}})\\
&\ll_{l,n} & 
 (1+|\tau|)^A \int_\R ~(1+|y_n|)^{A+l}
e^{-\frac{\pi}{2}\left(|y_n|+|\tau -y_n|\right)}~dy_n
\ll_{l,n}  (1+|\tau|)^A \psi_1(A+l;\tau).
\end{eqnarray*}
We finish the recurrence and the proof 
of the lemma \ref{majoration} by using 
the case $n=1$ above. \cqfd \par

\medskip
Before stating the central lemma of this paper, we will 
first introduce some needed notations:\\
Let $a \in \R_+^*$, $q\in \N$, $\rho >0$, $\delta >0$. 
Let $I$ be a finite subset of $\R^n \setminus \{\zerob\}$ and 
$\u=(u(\alphab))_{\alphab \in I}$ be a finite sequence of positive integers.
Set for all $\delta',\rho'>0$, 
$${\cal D}(\delta'; \rho'):= \{ (s,\z) \in \C \times \C^n \mid 
\sigma >a-\delta' {\mbox{ et }} |\Re (z_i)| <\rho' ~\forall i=1,\dots,n
\}.$$
Let $L(s,\z)$ be a holomorphic function on ${\cal D}(2\delta;2 \rho)$.
Assume that there exist $l,l'>0$ such that we have 
uniformly in $(s;\z) \in {\cal D}(2\delta; 2 \rho)$:
$$ L(s;\z) \ll \prod_{i=1}^n (1+|y_i|)^l (1+|\tau
-y_1-\dots-y_n|)^{\sigma + l'}
e^{\frac{\pi}{2}\left(|\tau|-\sum_{i=1}^n |y_i|
-|\tau -\sum_{i=1}^n y_i|\right)}.$$
We denote by $I_0$ the set of  real numbers  which are the coordinates 
  of at least one element of $I$, and by 
$\Q(I_0)$ the number field generated by $I_0$.\\
We set also for all  $\rhob =(\rho_1,\dots,\rho_n) \in
]-\rho,+\rho[^n$ 
such that $\rho_1,\dots,\rho_n$ are $\Q(I_0)$-linearly independent 
and for all $\betab \in \R^n$, 
$$T_n(s)=T_n(a,q,I,\u,\rhob,\betab,L; s):=\frac{1}{(2\pi i)^n} 
\int_{\rho_1-i\infty}^{\rho_1+i\infty}\dots\int_{\rho_n-i\infty}^{\rho_n+i\infty}  
\frac{L(s;\z)~dz_1 \dots dz_n}{(s-a-\la \betab , \z \ra )^q
  ~\prod_{\alphab \in I} {\la \alphab, \z \ra}^{u(\alphab)}}$$
We define finally $\eps(\betab; L)$ by:\\
$\eps(\betab; L)=1$ if $\betab \in con^*(I)\setminus \{\zerob \}$
  and if there exist two analytic functions $H$ and $U$ on ${\cal D}(2\delta;2
  \rho)$  such that  $L(s;\z)= U(s;\z) H\left(s; \left(\la \alphab , \z \ra
    \right)_{\alphab \in I}\right)$ and  $H(s;\zerob)\equiv 0$.\\
Otherwise we set $\eps(\betab ;L)=0$.\\
We can now state the crucial lemma of this paper:  
\begin{Lemme}\label{fonda1}
\be 
\item 
$s\mapsto T_n(s)$ converges absolutely in 
$\{\sigma > a+\sum_{i=1}^n |\beta_i \rho_i|\}$;
\item 
There exists $\eta >0$ such that $s\mapsto T_n(s)$ has a meromorphic
continuation with moderate growth to the half-plane $\{\Re(s) > a-\eta
\}$ with at most a single pole at   $s=a$;
\item 
If $s=a$ is a pole of $T_n(s)$ then its order is at most $d_n,$ where \\
$d_n:= \sum_{\alphab \in I} u(\alphab) -rank(I) +q -\eps(\betab; L)$.
\ee
\end{Lemme}
{\bf {\large Proof of  lemma \ref{fonda1}:}}\\
To prove lemma \ref{fonda1}, we will proceed by induction on $n$.\\
Throughout the   discussion, we use the following notations. 
Given $\z = (z_1,\dots,z_n)$ we set:
 $\z':=(z_1,\dots,z_{n-1})$ and  $l(\z):=\frac{1}{z_n} \z'=\left(\frac{z_1}{z_n},\dots,
  \frac{z_{n-1}}{z_n}\right)$ if $z_n \neq 0$. \par 
{\bf Step 1: Study of the case $n=1$:}\\
Set $A=\prod_{\alpha \in I} \alpha^{u(\alpha)}$ 
and $c=\sum_{\alpha \in I} u(\alpha)$. We have:
$$T_1(s)
=\frac{1}{(2\pi i)} 
\int_{\rho_1-i\infty}^{\rho_1+i\infty} 
\frac{L(s;z)~dz}{(s-a-\beta z )^q
  ~\prod_{\alpha \in I} (\alpha z)^{u(\alpha)}}=
\frac{1}{(2\pi i)A} 
\int_{\rho_1-i\infty}^{\rho_1+i\infty} 
\frac{L(s;z)~dz}{(s-a-\beta z )^q z^c}.$$
From our assumptions and lemma \ref{majoration} it follows that 
$s\mapsto T_1(s)$ converges absolutely and defines a holomorphic
function with moderate growth in  
$\{\sigma > a-\eta\}$ where $\eta=\inf(-\beta \rho_1, \delta)$.\\ 
$\bullet$ If $\beta \rho_1 <0$ then $\eta >0$. This proves the lemma
in this case.\\
$\bullet$ if $\beta =0$ then  
$T_1(s)=\frac{1}{(2\pi i)A (s-a)^q} 
\int_{\rho_1-i\infty}^{\rho_1+i\infty} 
\frac{L(s;z)~dz}{z^c}$.
It follows that  $s\mapsto T_1(s)$ 
has a meromorphic continuation with moderate growth to the half-plane 
$\{\sigma > a-\delta\}$ with at most one pole at  $s=a$ of order
at most $q\leq d_1$.\\
$\bullet$ We assume $\beta \rho_1 >0$:\\ 
The residue theorem and lemma \ref{majoration} imply that 
for $\sigma > a+\beta \rho_1$:\\ 
$T_1(s)=T_1'(s)+T_1''(s)$
where $T_1'(s)=\frac{1}{(2\pi i)A} 
\int_{-\rho_1-i\infty}^{-\rho_1+i\infty} 
\frac{L(s;z)~dz}{(s-a-\beta z )^q z^c}$ and $T_1''(s)=\frac{1}{A} 
Res_{z=0} \frac{L(s;z)}{(s-a-\beta z)^q z^c}$.\\
It is clear that $s\mapsto T_1'(s)$ converges absolutely  and
defines a meromorphic function with moderate growth in the half-plane
$\{\sigma > a-\eta\}$ where $\eta =\inf (\delta , \beta \rho_1)$.\\ 
If $c=0$ then $T_1''(s)\equiv 0$. Thus $T_1(s)=T_1'(s)$ satisfies  the
conclusions of  lemma \ref{fonda1}.\\
We assume now that $c\geq 1$. An easy computation shows that  
$$T_1''(s)=\frac{1}{A} \sum_{k=0}^{c-1} \frac{{-q \choose k} (-\beta)^k}{(c-1-k)!}~
\frac{{\partial_z}^{(c-1-k)}L(s;0)}{(s-a)^{q+k}}.$$
We deduce that $T_1(s)$ has a meromorphic continuation with moderate
growth to $\{\sigma > a-\eta\}$ with at most one pole at $s=a$
of order at most: 
\be
\item
ord$_{s=a} T_1(s) \leq q+c-1 =q+\sum_{\beta \in I} u(\beta) -rank(I)
=d_1$ \qquad \ if $L(s,0)\neq 0$ 
\newline (because in this case $\eps (\betab;L)=0$);
\item 
ord$_{s=a} T_1(s) \leq q+c-2 =q+\sum_{\beta \in I} u(\beta) -rank(I)-1
\leq  d_1$ \ \ if $L(s,0)= 0$;
\ee
Therefore the lemma \ref{fonda1} is true for $n=1$.\par
{\bf Step 2: Let $n \geq 2$. We assume that lemma \ref{fonda1}
 is true for any $\R^k, k \le n-1$. We will show that it 
also remains   true for $\R^n$:}\\
{\bf We  justify this assertion  by  induction 
on the integer} \\ 
$$h=h(\betab,\rhob):=\#\left\{i\in \{1,\dots,n\} \mid \beta_i \rho_i \geq
  0\right\} \in \{0,\dots,n\}.$$
$\bullet$ {\bf Proof of  lemma \ref{fonda1} for $h=0$:}\\
Since $h=0$ then for each $i=1,\dots,n,$ $\beta_i \rho_i <0$. 
It follows from lemma \ref{majoration} that 
$s\mapsto T_n(s)$ converges absolutely  and defines a holomorphic
function with moderate growth in the half-plane $\{\sigma >a-\eta\}$ 
where $\eta =\inf (-\la \betab,\rhob \ra, \delta) >0$. 
Thus lemma \ref{fonda1}  is also true in this case.\par
$\bullet$ {\bf Let $h\in \{1,\dots,n\}$. We assume that 
 lemma \ref{fonda1} is true for $h(\betab,\rhob)\leq h-1$. 
We will prove that it remains true for $h(\betab,\rhob)=h$ :}\\
If $\betab  =\zerob$ then 
$T_n(s)=\frac{1}{(2\pi i)^n (s-a)^q} 
\int_{\rho_1-i\infty}^{\rho_1+i\infty}\dots\int_{\rho_n-i\infty}^{\rho_n+i\infty}  
\frac{L(s;\z)~dz_1 \dots dz_n}{\prod_{\alphab \in I} {\la \alphab, \z
    \ra}^{u(\alphab)}}$. Since the $\rho_i$ are linearly 
independent   over  $\Q(I_0),$  lemma \ref{majoration} 
implies that lemma \ref{fonda1} is true in this case, in the sense that there is at most  
one  pole at $s=a$ of order at most $q\leq d_n$.\\
If $\betab \neq \zerob$ and $\beta_i \rho_i \leq 0$ for all 
$i=1,\dots,n$, then there exists $i_0$ such that  
$\beta_{i_0} \rho_{i_0} <0$. In this case, it is also easy to see  that 
$s\mapsto T_n(s)$ is  {\it holomorphic} with moderate growth in  
$\{ \sigma >a -\eta\}$ where $\eta =\inf (\delta , -\beta_{i_0}
\rho_{i_0}) >0$. It follows that lemma \ref{fonda1} is also true in
this case.\\
So to finish the proof of lemma \ref{fonda1} it suffices 
to consider the case where there exists $i\in \{1,
\dots,n\}$ such that $ \beta_i \rho_i >0$.\\
Without loss of generality we can assume that $\beta_n \rho_n >0$.\\
Set $\ds J:=\left\{\alphab \in I \mid \alpha_n \neq 0 {\mbox { and  }} 
\left|\frac{\alpha_1 \rho_1+\dots+\alpha_{n-1}
    \rho_{n-1}}{\alpha_n}\right|<|\rho_n|\right\}.$\\
Consider the equivalence relation ${\cal R}$ defined on $J$ by:
$\alphab {\cal R} \gammab$ iff $ \alpha_n \gammab =\gamma_n 
\alphab$. Denote by $J_1,\dots, J_r$ its equivalence classes  (they form a partition of $J$).\\
Choose for each $k=1,\dots,r,$ an element $\alphab^k \in J_k$ and 
set $c_k:=\sum_{\alphab \in J_k} u(\alphab)$.\\
Since $\rho_1,\dots,\rho_n$ are $\Q(I_0)$-linearly 
independent, it follows from 
the residue theorem and   lemma \ref{majoration}
that there exist constants $A_1,\dots, A_r
\in \R$ such that:
\begin{equation}\label{formuleresidu}
\forall \sigma > a+\sum_{i=1}^n |\beta_i \rho_i|, 
\quad T_n(s)=T_n^0(s)+\sum_{k=1}^r A_k T_{n-1}^k(s)
\end{equation}
where 
$$T_n^0(s):=\frac{1}{(2\pi i)^n} 
\int_{\rho_1-i\infty}^{\rho_1+i\infty}\dots\int_{\rho_{n-1}-i\infty}^{\rho_{n-1}+i\infty}
\int_{-\rho_n-i\infty}^{-\rho_n+i\infty}  
\frac{L(s;\z)~dz_1 \dots dz_n}{ 
(s-a-\la \betab , \z\ra)^q \prod_{\alphab \in I}{\la \alphab, \z
    \ra}^{u(\alphab)}}$$
and for each $k=1,\dots,r$ :
\begin{eqnarray*}
T_{n-1}^k(s)&:=&\frac{1}{(2\pi i)^{n-1}} 
\int_{\rho_1-i\infty}^{\rho_1+i\infty}\dots\int_{\rho_{n-1}-i\infty}^{\rho_{n-1}+i\infty}
\frac{1}{(c_k-1)!} \\
& & \left(\frac{\partial}{\partial z_n}\right)^{c_k-1}
\frac{L(s;\z)}{(s-a-\la \betab , \z\ra)^q 
\prod_{\alphab \in I\setminus J_k}{\la \alphab, \z \ra}^{u(\alphab)}}
\Big|_{z_n=-\la l(\alphab^k), \z'\ra}~dz_1 \dots dz_{n-1}.
\end{eqnarray*}
Since $h\left(\betab , (\rho_1,\dots,\rho_{n-1},-\rho_n)\right)
=h(\betab ,\rhob) -1 =h-1$, the induction hypothesis for $h-1$ 
implies that $s\mapsto T_n^0(s)$ satisfies the conclusions of 
lemma \ref{fonda1}. 
So to conclude, it is enough to prove lemma \ref{fonda1} for each 
$s\mapsto T_{n-1}^k(s)$. 

We then choose and fix any $k\in \{1,\dots,r\}$  for the rest of the discussion.\\
An easy computation shows that:
\begin{equation}\label{formuleresidubis}
T_{n-1}^k(s)=\sum_{u+v+\sum_{\alphab \in I\setminus J_k} k_{\alphab}=c_k -1} 
w\left(u,v,(k_{\alphab})\right) R_k\left(u,v,(k_{\alphab});s\right)
\end{equation}
where   $u,v$ and the $k_{\alphab}$   are in 
  $\N_0 $,  each  $w\left(u,v,(k_{\alphab})\right)\in \R$ and  
\begin{eqnarray}\label{rks}
R_k\left(u,v,(k_{\alphab});s\right)&:=& 
\frac{1}{(2\pi i)^{n-1}} 
\int_{\rho_1-i\infty}^{\rho_1+i\infty}\dots\int_{\rho_{n-1}-i\infty}^{\rho_{n-1}+i\infty}\\
& & \frac{ \frac{\partial^u L}{\partial {z_n}^u} (s;\z',-\la
  l(\alphab^k) , \z'\ra )~dz_1 \dots dz_{n-1}}{
\left(s-a-\la \betab' - \beta_n l(\alphab^k),
\z'\ra\right)^{q+v}
 \prod_{\alphab \in I\setminus J_k} {\la \alphab' - \alpha_n l(\alphab^k), \z'
    \ra}^{u(\alphab)+k_{\alphab}}}. \nonumber
\end{eqnarray}
So to conclude it suffices to prove lemma \ref{fonda1} for each 
$R_k\left(u,v,(k_{\alphab});s\right)$.\\
We fix now  $u, v \in \N_0 $ and $\left(k_{\alphab}\right)_{\alphab \in
  I\setminus J_k}$ such that 
$\ds u+v+\sum_{\alphab \in I\setminus J_k} k_{\alphab}=c_k -1$.\\
The induction hypothesis   implies that there exists $\eta
>0$ such that $s\mapsto R_k\left(u,v,(k_{\alphab});s\right)$ has
meromorphic continuation with moderate growth to the half-plane 
$\{\sigma  > a-\eta\}$ with at most  one   pole at $s=a$ of order 
at most 
\begin{equation}\label{majordre}
ord_{s=a} R_k\left(u,v,(k_{\alphab});s\right) \leq \bigg(\sum_{\alphab \in I \setminus J_k }
  u(\alphab)+k_{\alphab}\bigg) -rank(V) +(q+v) -
\eps\left(\betab' - \beta_n l(\alphab^k); \tilde{L}_u\right),
\end{equation}
where \ \  
$V:=\{\alphab' - \alpha_n l(\alphab^k) \mid \alphab
\in I \setminus J_k\}$ and 
$\tilde{L}_u (s;\z'):=\frac{\partial^u L}{\partial {z_n}^u} 
\left(s;\z',-\la l(\alphab^k), \z'\ra \right)$.\\
Set  $\tilde {V}:= \left\{\alphab -\frac{\alpha_n}{\alpha_n^k}
{\alphab^k}\mid \alphab \in I \setminus J_k\right\}
=\{\left(\alphab' - \alpha_n l(\alphab^k),
  0\right) \mid \alphab \in I \setminus J_k\}$.\\ 
It is clear that $rank(\tilde{V})=rank(V)$. 
Moreover it follows from the definition of $\alphab^k$ that $\alpha_n^k \neq 0,$
and therefore  
$\alphab^k \not \in Vect_\R(\tilde{V})$. We deduce that:
\begin{equation} \label{majrang} 
rank(I)=rank(\tilde{V}\cup \{\alphab^k\})=
rank(\tilde{V})+1=rank(V)+1.
\end{equation}
So it follows from (\ref{majordre}) and (\ref{majrang}) that: 
\begin{eqnarray*}
ord_{s=a} R_k\left(u,v,(k_{\alphab});s\right)  
&\leq & 
\big(\sum_{\alphab \in I \setminus J_k } u(\alphab)
+k_{\alphab}\big) -rank(V) +q+v-\eps\left(\betab' - \beta_n l(\alphab^k); \tilde{L}_u\right)\\
&= & c_k -1+\sum_{\alphab \in I \setminus J_k } u(\alphab) -rank(V)+q -u
-\eps\left(\betab' - \beta_n l(\alphab^k); \tilde{L}_u\right)\\
&\leq & \sum_{\alphab \in I} u(\alphab)-\left( rank(V)+1\right) +q 
-\eps\left(\betab' - \beta_n l(\alphab^k); \tilde{L}_u\right) -u\\
&\leq & \sum_{\alphab \in I} u(\alphab)-rank(I) +q 
-\eps\left(\betab' - \beta_n l(\alphab^k); \tilde{L}_u\right)
-u.
\end{eqnarray*}
Thus from the definition of $d_n$ we see that:
\begin{equation}\label{majrangbis}
ord_{s=a}  R_k\left(u,v,(k_{\alphab});s\right) \leq  d_n 
+\eps(\betab;L)-\eps\left(\betab' - \beta_n l(\alphab^k); \tilde{L}_u\right) -u.
\end{equation}
If $\eps(\betab;L)-\eps\left(\betab' - \beta_n l(\alphab^k); \tilde{L}_u\right) \leq 0$
or $u \neq 0$, 
then $\ds ord_{s=a} R_k\left(u,v,(k_{\alphab});s\right)  \leq  d_n$, 
which completes the proof.\\ 
So, we now assume   that $u=0$ and 
$\eps(\betab;L)-\eps\left(\betab' - \beta_n l(\alphab^k); \tilde{L}_0\right) >0$, that is,\\ 
$u=0, \quad \eps(\betab;L)=1 {\mbox { and }}  
\eps\left(\betab' - \beta_n l(\alphab^k); \tilde{L}_0\right)=0
$.\\
Therefore we assume that $\betab \in con^* (I) \setminus \{\bold 0\}$, and $L$ is of the form  
$L(s;\z)= U(s;\z) H\left(s; \left(\la \alphab , \z \ra
    \right)_{\alphab \in I}\right)$ with $H(s;\zerob)\equiv 0$.\\ 
It is then clear that ${\tilde L}_0$ is also of the form:
\begin{equation}\label{sauvetout} 
{\tilde L}_0(s,\z')= {\tilde U}(s;\z') {\tilde H}\left(s; \left(\la \mub , \z' \ra
    \right)_{\mub \in V}\right) {\mbox { with }} {\tilde H}(s;\zerob)\equiv 0.
\end{equation} 
This and the fact that 
$ \eps\left(\betab' - \beta_n l(\alphab^k); \tilde{L}_0\right)=0$
imply that we have necessarily: 
\begin{equation}\label{con*}
  \betab'-\frac{\beta_n}{\alpha_n^k} {(\alphab^k)}' 
\not \in con^*(V)\setminus \{\zerob\}.
\end{equation} 
Since $\betab \in con^* (I) \setminus \{\bold 0\}$, there exists  
$  \{\lambda_{\alphab}\}_{\alphab \in I} \subset \R_+^*$   such
that $\betab = \sum_{\alphab \in I} \lambda_{\alphab} \alphab$.
This implies  
$\betab'= \sum_{\alphab \in I} \lambda_{\alphab} \alphab'$ and  
$\beta_n = \sum_{\alphab \in I} \lambda_{\alphab} \alpha_n$. 
We deduce that:
\begin{eqnarray*}
\betab'-\beta_n l(\alphab^k)
&=& \sum_{\alphab \in I} \lambda_{\alphab} \alphab' 
- \sum_{\alphab \in I} \lambda_{\alphab}
\alpha_n l(\alphab^k)
= \sum_{\alphab \in I} \lambda_{\alphab} 
\left(\alphab' - \alpha_n l(\alphab^k)\right)\\
&=& \sum_{\alphab \in I\setminus J_k} 
\lambda_{\alphab} 
\left(\alphab' - \alpha_n l(\alphab^k)\right)
~({\mbox {because }} \alphab' - \alpha_n l(\alphab^k) =0 {\mbox
  { if }} \alphab \in J_k ).
\end{eqnarray*}
It follows that $\betab'-\beta_n l(\alphab^k)
\in con^*(V)$. Therefore (\ref{con*}) gives 
$\betab'-\beta_n l(\alphab^k) 
=\zerob$. \\
It follows then from (\ref{rks})  that 
\begin{equation} \label{majord+}
ord_{s=a} R_k\left(0,v,(k_{\alphab});s\right) \leq q+v =q+c_k -1-\sum_{\alphab \in I\setminus
  J_k} k_{\alphab} 
\leq  q +\sum_{\alphab \in J_k} u(\alphab) -1.
\end{equation}
But $\betab'-\beta_n l(\alphab^k)
=\zerob$ implies also that 
$\betab -\frac{\beta_n}{\alpha_n^k} {\alphab^k} =
\left(\betab'-\beta_n l(\alphab^k), 0\right)
=\zerob$. Therefore we get $\ds \betab =\frac{\beta_n}{\alpha_n^k}
{\alphab^k}$.
In particular, since $\betab \neq \zerob$, we have $\beta_n \neq
0$ and thus $Vect_\R(J_k)=\R \alphab^k =\R \betab$.
Since $\betab \in con^*(I)$, it follows that 
$rank (I) = rank (I \setminus \{\betab\})$, which, we therefore now see, implies     
$rank (I) = rank(I\setminus J_k),$ unless $I=J_k$. \\
{\bf Assume first that $I\neq J_k$:} \\
In this event, combining (\ref{majord+}) with the fact 
that $rank(I\setminus J_k)=rank(I)$, we conclude:
\begin{eqnarray*}
ord_{s=a} R_k\left(0,v,(k_{\alphab});s\right) 
&\leq & 
q+\sum_{\alphab \in I} u(\alphab) -rank(I)-
\sum_{\alphab \in I\setminus J_k} u(\alphab)+rank(I) -1\\
& \leq & q+\sum_{\alphab \in I} u(\alphab) -rank(I)-1
-\left(\sum_{\alphab \in I\setminus J_k} u(\alphab)-
rank(I\setminus J_k)\right) \\
& \leq & q+\sum_{\alphab \in I} u(\alphab) -rank(I)-1
-\left(\#\left(I\setminus J_k\right) -rank(I\setminus J_k)\right)\\
&\leq & q+\sum_{\alphab \in I} u(\alphab) -rank(I)-1 \leq d_n.
\end{eqnarray*}
{\bf Assume that $I=J_k$:}\\
It is then clear that $V=\emptyset$ and it follows from (\ref{sauvetout})
that ${\tilde L}_0(s;\z')\equiv 0$.
Since $u=0$,  (\ref{rks}) implies that 
$R_k\left(0,v,(k_{\alphab});s\right)\equiv 0$,   in which case,  it is
obvious that we have  
$ord_{s=a} R_k\left(0,v,(k_{\alphab});s\right) \leq  d_n$.\par 
We conclude that for any $u, v, (k_{\alphab}),\ $    
$ord_{s=a} R_k\left(u,v,(k_{\alphab});s\right) \leq  d_n$. 
This finishes the induction argument on $h,$  therefore, also on $n$, and   
completes the proof of lemma \ref{fonda1}.\cqfd
\begin{Lemme}\label{fonda3}
Let $\a =(a_1,\dots,a_r) \in \R_+^{*r}$ and  $a=|\a|=a_1+\dots+a_r$.
Let $I$ be a finite nonempty subset of $\R_+^r \setminus\{\zerob\}$, 
$\u=\left(u(\betab)\right)_{\betab \in I}$ a finite sequence of
positive integers and $\h=(h_1,\dots, h_r) \in \R_+^{*r}$.
Assume that: $\unb \in con (I)$ and that $\la \alphab , \a \ra =1$ for
all $\alphab \in I$.\\
Let  $\rhob \in \R_+^{*r}$. For $\sigma =\Re(s)> a+|\rhob|$ set:
$$ {\cal R}(s):=
\frac{1}{(2\pi i)^r} 
\int_{\rho_1-i\infty}^{\rho_1+i\infty}\dots
\int_{\rho_r-i\infty}^{\rho_r+i\infty}  
\frac{\Gamma(s-a-z_1-\dots-z_r)~\prod_{i=1}^r \Gamma (a_i+z_i) ~
d\z}
{\Gamma(s)~\prod_{k=1}^r h_k^{a_k+z_k}~\prod_{\alphab \in I} 
{\la \alphab, \z \ra}^{u(\alphab)}}.$$
Then there exists $\eta >0$ such that 
$s\mapsto {\cal R}(s)$ has a meromorphic 
continuation to the half-plane $\{\sigma >a-\eta\}$ 
with exactly  one pole at $s=a$ of order  
$\rho_0:= \sum_{\betab \in I}
u(\betab) - rank\left(I \right) +1$. Moreover we have 
$\ds {\cal R}(s)\sim_{s \rightarrow a} 
\frac{A_0(I;\u;\h)}{(s-a)^{\rho_0}}\,,$ where $A_0(I;\u;\h )$ is the
volume constant (see \S 2.3) associated to $I$, $\u$ and $\h$.
\end{Lemme}
{\bf Proof of lemma \ref{fonda3}:}\\
Set $q:=\sum_{\alphab \in I} u(\alphab)$.  
We define $\alphab^1,\dots,\alphab^q \in \R_+^r \setminus \{\zerob\}$
by: \\
$\{\alphab^i \mid i=1,\dots,q\}=I$ and $\forall \alphab \in I$ 
$\#\{i\in \{1,..,q\} \mid \alphab^i=\alphab \}=u(\alphab)$.\\
We define then  $\mub^1,\dots,\mub^r \in \R_+^q$ by: 
$\mu_i^k=\alpha_k^i$ $\forall i=1,\dots,n$ and 
$\forall k=1,\dots,r$.\\
Set  
$G(\X)=1+P_{(I;\u;\h)}(\X):=1+\sum_{k=1}^r h_k~\X^{\mub^k}$. 
$G$ is a generalized polynomial with positive coefficients.  
Moreover for any $i=1,\dots,q$ there exists $k\in \{1,\dots ,r\}$ such
that $\mu^k_i=\alpha^i_k \neq 0$. It follows that $P$ 
depends on all the variables $X_1,\dots,X_q$.\\
We will first prove that for $\sigma \gg 0$, ${\cal R}(s)=\int_{[1,+\infty[^n} G^{-s}(\x)~ d\x$.\\
Mellin's formula (\ref{mellinformula})
implies that for all $\sigma \gg 0$ :
\begin{eqnarray*}
{\cal R}(s)&:=&
\frac{1}{(2\pi i)^r} 
\int_{\rho_1-i\infty}^{\rho_1+i\infty}\dots
\int_{\rho_r-i\infty}^{\rho_r+i\infty}  
\frac{\Gamma(s-a-z_1-\dots-z_r)~\prod_{i=1}^r \Gamma (a_i+z_i) ~
d\z}
{\Gamma(s)~\prod_{i=1}^r h_i^{a_i+z_i}~\prod_{k=1}^q {\la \alphab^k, \z \ra}}\\
& = & 
\frac{1}{(2\pi i)^r} 
\int_{a_1+\rho_1-i\infty}^{a_1+\rho_1+i\infty}\dots
\int_{a_r+\rho_r-i\infty}^{a_r+\rho_r+i\infty}  
\frac{\Gamma(s-z_1-\dots-z_r)~\prod_{i=1}^r \Gamma (z_i) ~
d\z}
{\Gamma(s)~\prod_{i=1}^r h_i^{z_i}~\prod_{k=1}^q {\left(\la \alphab^k, \z \ra-1\right)}}\\
&=& 
\frac{1}{(2\pi i)^r} 
\int_{a_1+\rho_1-i\infty}^{a_1+\rho_1+i\infty}\dots
\int_{a_r+\rho_r-i\infty}^{a_r+\rho_r+i\infty}  
\Gamma(s-z_1-\dots-z_r)~{\Gamma(s)}^{-1}~\prod_{i=1}^r \Gamma (z_i)~ 
\prod_{i=1}^r h_i^{-z_i}\\
& & \qquad \qquad \times 
\left(\int_{ [1,+\infty[^q} {\prod_{k=1}^q x_k^{-\la \alphab^k, \z
      \ra}}~dx_1\dots dx_q\right)
dz_1 \dots dz_r\\
&=& 
\frac{1}{(2\pi i)^r} 
\int_{a_1+\rho_1-i\infty}^{a_1+\rho_1+i\infty}\dots
\int_{a_r+\rho_r-i\infty}^{a_r+\rho_r+i\infty}  
\Gamma(s-z_1-\dots-z_r)~{\Gamma(s)}^{-1}~\prod_{i=1}^r \Gamma (z_i)~
\prod_{i=1}^r h_i^{-z_i}\\
& & \qquad \qquad \times 
\left(\int_{ [1,+\infty[^q} {\prod_{i=1}^r {\x}^{-z_i
      \mub^i}}~dx_1\dots dx_q\right)
dz_1 \dots dz_r\\
& = &
\int_{ [1,+\infty[^q}\Big[
\frac{1}{(2\pi i)^r} 
\int_{a_1+\rho_1-i\infty}^{a_1+\rho_1+i\infty}\dots
\int_{a_r+\rho_r-i\infty}^{a_r+\rho_r+i\infty}  
\Gamma(s-z_1-\dots-z_r)~{\Gamma(s)}^{-1}~\prod_{i=1}^r \Gamma (z_i) \\
& & \qquad \qquad \times {\prod_{i=1}^r {\left(h_i~{\x}^{\mub^i}\right)}^{-z_i}}
dz_1 \dots dz_r \Big]dx_1\dots dx_q \\
&= & \int_{ [1,+\infty[^q}{\bigg(1+\sum_{i=1}^r
    h_i~{\x}^{\mub^i}\bigg)}^{-s}~dx_1\dots dx_q 
=\int_{[1,+\infty[^n} G^{-s}(\x)~ d\x.
\end{eqnarray*}
So to conclude it suffices to check that $s\mapsto
Y(G;s):=\int_{[1,+\infty[^n} G^{-s}(\x)~ d\x$ satisfies   
the conclusions of lemma \ref{fonda3}.\\
Let ${\cal E}^{\infty}(G)=conv \left(supp(G)-\R_+^q\right)$ denote 
the Newton polyhedron at infinity of  $G$.
Denote  by $G_0$ the smallest face that meets the diagonal. 
It follows from  Sargos' result (see \S 2.3) that there exists  
$\eta >0$ such that $Y(G;s)$ has a meromorphic continuation to the
half-plane   $\{\sigma>\sigma_0 -\eta\}$ (where $\sigma_0 = \sigma_0(P)$)
with moderate growth and exactly one pole at $s=\sigma_0$ of order 
$\rho_0:=codim G_0$. Moreover
$\sigma_0$ is characterized geometrically by: 
${\sigma_0}^{-1} \unb =\Delta \cap G_0 =
\Delta \cap {\cal E}^{\infty}(G)$ and  
$Y(G;s) \sim_{s \rightarrow \sigma_0} 
\frac{A_0(G)}{(s-\sigma_0)^{\rho_0}}$ where  
$A_0(G)$ is the Sargos constant associated to the polynomial $G$. 
It is easy to see that in our case $A_0(G)$ is equal 
to the volume constant $A_0(I;\u;\h)$ associated to $I$, $\u$ and $\h$.\\
By our hypothesis,  we have 
$\unb \in con (I)=
con \left(\{ \alphab^k \mid k=1,\dots,q\}\right)$. Thus there exists  
$\cb=(c_1,\dots,c_q) \in \R_+^q \setminus \{\zerob\}$ such that 
$\ds \unb =\sum_{k=1}^q c_k \alphab^k$. 
It follows that: \\
$$\forall i=1,\dots ,r \quad \la \cb,\mub^i \ra =
\sum_{k=1}^q c_k \mu_k^i=\sum_{k=1}^q c_k \alpha_i^k =1.$$
Since   each  $\mub^i \in supp (G),$ we conclude that   
$L_{\cb}:=\{ \x \in \R^q \mid \la \cb , \x \ra =1\}$ is a support  plane 
of ${\cal E}^{\infty}(G).$ 
Thus: 
\begin{equation}\label{gammaface}
F_{\cb}^\infty :=L_{\cb} \cap {\cal E}^{\infty}(G)=
conv\left(\{\mub^i \mid i=1,\dots,r\}\right) 
{\mbox { is a face of the polyhedron }} {\cal E}^{\infty}(G).
\end{equation}
By our hypothesis we know that   
$\sum_{k=1}^r a_k \mu_i^k= \la \a ,\alphab^i \ra =1 \ \forall i=1,\dots ,q,$  which implies  
 $\ds  \frac{1}{a}\unb =\sum_{k=1}^r \frac{a_k}{a} \mub^k \in conv\left(\{\mub^i \mid i=1,\dots,r\}\right).$  Thus,  $\frac 1a
\unb \in F_{\cb}^\infty \cap
\Delta,$ that is, the face $F_{\cb}^\infty$ must meet the diagonal at 
$\frac{1}{a} \unb$. It follows that $\sigma_0=a$ and 
that $G_0 \subset F_{\cb}^\infty$. Hence we deduce that: 
\begin{eqnarray*}
ord_{s=a}Y(G;s)&=& codim G_0 \geq  codim F_{\cb}^\infty =q-dim F_{\cb}^\infty \\
&\geq & q -rank (\left(\{\mub^i \mid i=1,\dots,r\}\right)+1 
= q-rank (\left(\{\alphab^i \mid i=1,\dots,q\}\right)+1 \\
&\geq & q-rank (I) +1= \sum_{\alphab \in I} u(\alphab) -rank(I)+1.
\end{eqnarray*}
But lemma \ref{fonda1} implies that 
$ord_{s=a}Y(G;s)=ord_{s=a}{\cal R}(s) 
\leq 1+\sum_{\alphab \in I} u(\alphab) -rank(I)$.  So we 
have   
$ord_{s=a} Y(G;s)=\sum_{\alphab \in I} u(\alphab) -rank(I)+1$. 
This completes the proof of lemma \ref{fonda3}. \cqfd
\subsection{Proof of theorem \ref{analytic}}
Let $\cb\in \Sigma_f (P)=\Sigma_f  \cap con^*\left(supp(P)\right)$ and 
$\T=(I_{\cb};\u )$ the polar type of $f$ in it.\\
Since  $\sum \frac{|f(m_1,\dots,m_n)|}{(m_1\dots m_n)^t}<+\infty$ for  
$t=1+\sup_{i} c_i$, we certainly have   
\begin{equation}\label{majf}
f(m_1,\dots,m_n)\ll (m_1\dots m_n)^t {\mbox { uniformly in }} \m \in
\N^{n}. 
\end{equation}
Let $P\in \R_+[X_1,\dots,X_n]$ be a homogeneous polynomial of degree 
$d\geq 1$ which depends on all the variables  $X_1,\dots,X_n$. Set 
$P(\X):=\sum_{k=1}^r b_k \X^{\gammab^k}$ ($b_k \in \R_+^*$ $\forall k$).\\
Since $\cb \in con^*\left(supp(P)\right),$   there exists  
$\a =(a_1,\dots ,a_r)\in \R_+^{*r}$, such that $\cb =\sum_{k=1}^r a_k
\gammab^k$. It follows that we have uniformly in
$\x \in [1,+\infty[^n$:
\begin{equation}\label{majp+}
\left(x_1^{c_1}\dots x_n^{c_n}\right)^{1/|\a|}
=\left(\x^{\sum_{k=1}^r (a_k/|\a|)~ \gammab^k }\right)\ll 
\sum_{k=1}^r \x^{\gammab^k} \ll
P(\x) \ll \left(x_1\dots x_n\right)^d.
\end{equation}
It is clear that (\ref{majf}) and (\ref{majp+}) imply that the series  
$Z(f;P;s):=\sum_{\m \in \N^{n}}
\frac{f(m_1,\dots,m_n)}{P(m_1,\dots,m_n)^{s/d}}$ 
has an abscissa of convergence $\sigma_0 <+\infty$.\\
Moreover, (\ref{majp+}) also implies that $\forall M\in \N$ 
we have uniformly in $\x \in [1,+\infty[^{n}$
and $\s \in \C$:\\
$
P(\x)^{-s/d}
= \sum_{k=0}^M (-1)^k {-s/d \choose
  k}~{\left(1+P(\x)\right)}^{-(s+dk)/d}
+ 
{\cal O}\left((1+|s|^{M+1}) \left(1+P(\x)\right)^{-(\Re(s)+d
    M+d)/d}\right).
$\\
It follows that for $M\in \N$ and $\sigma >\sigma_0$:
\begin{equation}\label{ZP+1}
Z(f;P;s)
= \sum_{k=0}^M (-1)^k {-s/d \choose
  k}~Z(f;1+P;s+dk)
+ 
{\cal O}\left((1+|s|^{M+1}) Z(|f|;1+P; \sigma +d
    M+d)\right).
\end{equation} 
Thus,  it suffices to prove the assertion of the theorem for  
$Z(f;1+P;s)$. \par
Let $\alphab^1,\dots, \alphab^n$
be $n$ elements of $\N_0^r\setminus \{\zerob\}$ defined by:
$\alpha^i_k=\gamma^k_i$  for all $i=1,\dots,n $ and    
$k=1,\dots,r $.
It follows from the assumptions on $f$ that there exists $\eps_1  >1$ such that 
\begin{equation}\label{Tz} 
\z \mapsto T(\T;P;\z):=H_{\cb}\left(f; \la \alphab^1,\z \ra ,\dots , \la
  \alphab^n,\z \ra\right)
\end{equation}
has a holomorphic continuation with moderate growth to 
$\{\z \in \C^r \mid \forall i ~\Re (z_i) >-\eps_1\}$.\\
By using   Mellin's formula (\ref{mellinformula}) we obtain that for
any $\rhob \in \R_+^{*n}$ and $\sigma > \sup \left(\sigma_0, d |\rhob|\right)$:
\begin{eqnarray*}
& & \Gamma (s/d) ~Z(f;1+P;s)
= 
\Gamma (s/d) ~\sum_{\m \in \N^{n}} \frac{f(m_1,\dots ,
  m_n)}{\left(1+\sum_{k=1}^r b_k \m^{\gammab^k}\right)^{s/d}}\\
&=& 
\sum_{\m \in \N^{n}}
\frac{1}{(2\pi i)^r} 
\int_{\rho_1-i\infty}^{\rho_1+i\infty}\dots
\int_{\rho_r-i\infty}^{\rho_r+i\infty}  
\Gamma(s/d -z_1-\dots-z_r)~\prod_{i=1}^r \frac{\Gamma (z_i)}{b_i^{z_i}} ~
\frac{f(m_1,\dots ,
  m_n)}{\left(\prod_{k=1}^r \m^{z_k \gammab^k}\right)} d\z\\
& = & 
\sum_{\m \in \N^{n}}
\frac{1}{(2\pi i)^r} 
\int_{\rho_1-i\infty}^{\rho_1+i\infty}\dots
\int_{\rho_r-i\infty}^{\rho_r+i\infty}  
\Gamma(s/d -|\z|)~\prod_{i=1}^r \frac{\Gamma (z_i)}{b_i^{z_i}} ~
\frac{f(m_1,\dots ,
  m_n)}{\left(\prod_{i=1}^n m_i^{\la \alphab^i, \z \ra }\right)}
d\z \\
&=& 
\frac{1}{(2\pi i)^r} 
\int_{\rho_1-i\infty}^{\rho_1+i\infty}\dots
\int_{\rho_r-i\infty}^{\rho_r+i\infty}  
\Gamma(s/d -|\z|)~\prod_{i=1}^r \frac{\Gamma (z_i)}{b_i^{z_i}} ~
{\cal M}\left(f; \la \alphab^1, \z \ra , \dots, \la \alphab^n, \z \ra \right)~
d\z.
\end{eqnarray*}
For all $\betab \in I_{\cb}$, we set:\\
 $\mu(\betab):=\sum_{i=1}^n
\beta_i \alphab^i,$ \ \   $I_{\cb}^*=\{\mu(\betab) \mid \betab \in I_{\cb}\},$ \ \ and  
$\u^*=\left(u^*(\etab)\right)_{\etab \in I_{\cb}^*},$\\ where \\
$u^*(\etab)=\sum_{\{\betab \in I_{\cb}; ~\mu(\betab)=\etab\}}
u(\betab)$ $\forall \etab \in I_{\cb}^*$.\\
The relation $\cb =\sum_{k=1}^r a_k \gammab^k$ implies that 
$\la \alphab^i , \a \ra =c_i ~\forall i$. 
Therefore it follows from the previous computations that for all 
$\rhob \in \R_+^{*n}$ and for all 
$\sigma > \sup \left(\sigma_0, d (|\a|+|\rhob|)\right)$,
\begin{eqnarray*}
& & \Gamma (s/d) ~Z(f;1+P;s)\\
&=& 
\frac{1}{(2\pi i)^r} 
\int_{a_1+\rho_1-i\infty}^{a_1+\rho_1+i\infty}\dots
\int_{a_r+\rho_r-i\infty}^{a_r+\rho_r+i\infty}  
\Gamma(s/d -|\z|)~\prod_{i=1}^r \frac{\Gamma (z_i)}{b_i^{z_i}} ~
{\cal M}\left(f; \la \alphab^1, \z \ra , \dots, \la \alphab^n, \z \ra
  \right) 
~
d\z \\
&=&
\frac{1}{(2\pi i)^r} 
\int_{\rho_1-i\infty}^{\rho_1+i\infty}\dots
\int_{\rho_r-i\infty}^{\rho_r+i\infty}  
\Gamma(s/d -|\a|-|\z|)~\prod_{i=1}^r \frac{\Gamma (a_i+z_i)}{b_i^{a_i+z_i}} ~
{\cal M}\left(f; c_1+\la \alphab^1, \z \ra , \dots, c_n+\la \alphab^n, \z \ra
  \right) d\z \\
&=& 
\frac{1}{(2\pi i)^r} 
\int_{\rho_1-i\infty}^{\rho_1+i\infty}\dots
\int_{\rho_r-i\infty}^{\rho_r+i\infty}  
\Gamma(s/d -|\a|-|\z|)~\prod_{i=1}^r \frac{\Gamma(a_i+z_i)}{b_i^{a_i+z_i}} ~
\frac{T(\T;P;\z)}{\prod_{\etab \in I_{\cb}^*} \left(\la\etab, \z
    \ra\right)^{u^*(\etab)}} ~d\z 
\end{eqnarray*}
We deduce from this that for any $\rhob \in \R_+^{*n}$ and  
$\sigma > \sup \left(\sigma_0, d (|\a|+|\rhob|)\right)$:
\begin{equation}\label{formuleclef}
Z(f;1+P;s)=
\frac{1}{(2\pi i)^r} 
\int_{\rho_1-i\infty}^{\rho_1+i\infty}\dots
\int_{\rho_r-i\infty}^{\rho_r+i\infty} 
\frac{ V(s;\z) ~d\z}{ \left(s/d -|\a|-\la \unb ,\z \ra \right) 
\prod_{\etab \in I_{\cb}^*} \left(\etab, \z
    \ra\right)^{u^*(\etab)}}
\end{equation}
where 
$$V(s;\z):=\left(s/d -|\a|-|\z|\right) 
\Gamma(s/d -|\a|-|\z|) {\Gamma(s/d )}^{-1}~
\left(\prod_{i=1}^r \Gamma(a_i+z_i) b_i^{-a_i-z_i}\right) ~
T(\T;P;\z).$$
By using the classical properties of the Euler $\Gamma$ 
function (in particular Stirling's formula) and (\ref{Tz}) 
it is easy to check that the assumptions of lemma {\ref{fonda1}} are satisfied. 
This   implies   that there exists $\eta >0$ such that 
$s\mapsto Z(f;1+P;s)$ has a meromorphic continuation with moderate 
growth to the half-plane $\{\sigma > d |\a|-\eta\}$ with at most one 
  pole at $s=d|\a|$.  Since  $|\cb|=\la \unb ,\cb \ra =\la \unb , \sum_{k=1}^r a_k \gammab^k\ra =
 \sum_{k=1}^r a_k |\gammab^k|= d |\a|$, 
the proof of the theorem follows immediately. \cqfd \par
{\bf Proof of   theorem \ref{partieprincipale}:}\\
We assume now that there exists 
$\cb \in \Sigma_f (P;\unb):={\cal F}(\Sigma_f)(\unb) 
\cap con^*\left(supp(P)\right)$ satisfying: 
\be
\item $\unb \in con^*(I_{\cb})$;
\item there exists a function $L$ such that: 
$\ds {\cal M}(f;\s)= L\left((\la \betab,\s\ra)_{\betab \in
 I_{\cb}}\right)$.
\item $H_{\cb}(f;\bold 0) \neq 0.$
\ee
We also continue to use the notations of the preceding paragraph.
Set 
$$U(s;\z):= \frac{\left(s/d -|\a|-|\z|\right) 
\Gamma(s/d -|\a|-|\z|)}{\Gamma(s/d ) }~\prod_{i=1}^r
\frac{\Gamma(a_i+z_i)}{b_i^{a_i+z_i}} {\mbox { and }} 
{\tilde T}(\z):= T(\T;P;\z)-T(\T;P;\zerob).$$
Then 
$\ds V(s;\z):= U(s;\z) ~T(\T;P;\z)= H_{\cb}(f;\zerob)~U(s;\z) + U(s;\z)~{\tilde
  T}(\z)$.\\
We note that the formula  (\ref{formuleclef})   can be written as follows: 
\begin{equation}\label{zdecomposition}
Z(f;1+P;s) = Z_1(s) +Z_2(s) \quad (\sigma \gg 0).
\end{equation}
where
\begin{eqnarray*}
Z_1(s)&=& \frac{H_{\cb}(f;\zerob)}{(2\pi i)^r} 
\int_{\rho_1-i\infty}^{\rho_1+i\infty}\dots
\int_{\rho_r-i\infty}^{\rho_r+i\infty} 
\frac{ U(s;\z) ~d\z}{\left(s/d -|\a|-\la \unb ,\z \ra \right) 
\prod_{\etab \in I_{\cb}^*} \left(\etab, \z
    \ra\right)^{u^*(\etab)}}\\
Z_2(s)&=& \frac{1}{(2\pi i)^r} 
\int_{\rho_1-i\infty}^{\rho_1+i\infty}\dots
\int_{\rho_r-i\infty}^{\rho_r+i\infty} 
\frac{ U(s;\z) {\tilde T}(\z)~d\z}{\left(s/d -|\a|-\la \unb ,\z \ra \right) 
\prod_{\etab \in I_{\cb}^*} \left(\etab, \z
    \ra\right)^{u^*(\etab)}}\,.
\end{eqnarray*}
$\bullet$ {\bf Study of  $Z_2(s)$:}\\
Since $\unb \in con^*(I_{\cb})$, there exists a set   
$\left\{t_{\betab}\right\}_{\betab \in I_{\cb}} \subset \R_+^*$ such that \\
$\unb =\sum_{\betab \in I_{\cb}} t_{\betab} \betab$ (i.e  
$\sum_{\betab \in I_{\cb}} t_{\betab} \beta_i =1$ $\forall i=1,\dots,n$). 
Consequently we have:
$$
\sum_{\betab \in I_{\cb}} t_{\betab} \mu(\betab)
= \sum_{\betab \in I_{\cb}} t_{\betab} \sum_{i=1}^n \beta_i
\alphab^i 
= 
\sum_{i=1}^n \big(\sum_{\betab \in I_{\cb}} t_{\betab}  \beta_i\big)
\alphab^i = \sum_{i=1}^n \alphab^i  
=\sum_{k=1}^r \big(\sum_{i=1}^n \gamma_i^k\big) \e_k
= d \unb.
$$
We conclude from this that $\unb \in con^*\left(I_{\cb}^*\right)\setminus
\{\zerob \}$.\par
Furthermore,  assumption 2 implies that there exists a function 
$L_1$ such that 
$\ds H_{\cb}(f;\s)= L_1\left((\la \betab,\s\ra)_{\betab \in
 I_{\cb}}\right)$. 
But for any $\betab \in I_{\cb}$ and for any $\z \in \C^r$ we have, 
$$\sum_{i=1}^n \beta_i \la \alphab^i,\z\ra =
\left\la \sum_{i=1}^n \beta_i \alphab^i,\z\right \ra= \la
\mu(\betab), \z\ra.$$
It follows that:
$T(\T;P;\z):=H_{\cb}\left(f; \la \alphab^1, \z\ra,\dots, \la
  \alphab^n, \z\ra\right)=
L_1\left(\left( \la \mu(\betab), \z\ra\right)_{\betab \in
    I_{\cb}}\right)$.
Consequently there exists a function ${\tilde L}$ such that:
${\tilde T}(\z)={\tilde L}\left((\la \etab,\z\ra)_{\etab \in
    I_{\cb}^*}\right)$. 
Since in addition we have $|\cb|=d|\a|$, ${\tilde T}(\zerob)=0$ and $\unb \in
con^*(I_{\cb}^*)$, it follows from lemma \ref{fonda1} that:
\begin{equation}\label{ordz2}
 ord_{s=|\cb|} Z_2(s)\leq \rho_0^* -1 {\mbox { where }}
 \rho_0^*:=\sum_{\etab \in I_{\cb}^*} u^*(\etab) -rank(I_{\cb}^*) +1 =\rho_0(\T_{\cb};P).
\end{equation} 
$\bullet$ {\bf Study of $Z_1(s)$ :}\\
It is easy to see that for $\sigma \gg 0$ :
$$Z_1(s)
=\frac{H_{\cb}(f;\zerob)}{(2\pi i)^r} 
\int_{\rho_1-i\infty}^{\rho_1+i\infty}\dots
\int_{\rho_r-i\infty}^{\rho_r+i\infty} 
\frac{ \Gamma(s/d -|\a|-|\z|) ~ \prod_{k=1}^r
\Gamma(a_k+z_k)~d\z}{\Gamma (s/d)~ \prod_{k=1}^r b_k^{a_k+z_k}  
\prod_{\etab \in I_{\cb}^*} \left(\la \etab, \z
    \ra\right)^{u^*(\etab)}}\,.$$
Moreover we know that 
$H_{\cb}(f;\zerob)\neq 0$, $\unb \in con(I_{\cb}^*),$ and  
$\la\etab ,\a\ra=\sum_{i=1}^n \beta_i \la\alphab^i, \a\ra= 
\sum_{i=1}^n \beta_i c_i=\la\betab , \cb\ra =1$ $\forall \etab =\mu (\betab) \in I_{\cb}^*$. 
Thus, it follows from  lemma \ref{fonda3} that $s=d|\a|=|\cb|$ is a
pole of $Z_1(s)$ of order  $\rho_0^*$ and that 
$\ds Z_1(s)\sim_{s\rightarrow |\cb|} \frac{H_{\cb}(f;\zerob) d^{\rho_0^*}
    A_0(I_{\cb}^*;\u^*;\b)}{\left(s-|\cb|\right)^{\rho_0^*}}\,,$  
where $A_0(I_{\cb}^*;\u^*;\b)>0$ is the {\it volume constant} 
associated to $I_{\cb}^*, \u^*$ and $\b$.
Combining this with  (\ref{ordz2}),
(\ref{zdecomposition}), and (\ref{ZP+1}), implies that 
$Z(f;P;s)$ has a pole at $s=|\cb|$ of order  
$\rho_0^*\,,$  and that \\
$Z(f;P;s)\sim_{s\rightarrow |\cb|} \frac{H_{\cb}(f; \zerob) d^{\rho_0^*}
    A_0(I_{\cb}^*;\u^*;\b)}{\left(s-|\cb|\right)^{\rho_0^*}}$. 
This completes the proof of   theorem \ref{partieprincipale}.\cqfd

\subsection{Proof of corollary \ref{mats+}}
Define a function $f:\N^n \rightarrow \C$ by:
$f(m_1,\dots,m_n)=\prod_{j=1}^n f_j(m_j)$. It is easy to see that 
$\ds {\cal M}(f;\s):=\sum_{m_1,\dots,m_n \geq 1} \frac{f(m_1,\dots,m_n)}{m_1^{s_1}\dots
  m_n^{s_n}}$  converges absolutely and satisfies\\ 
${\cal M}(f;\s):=\prod_{j=1}^n Z(f_j;s_j)$ in 
$\Omega:=\{\s \in \C^n \mid \sigma_i >c_i ~\forall i\}$. It follows that 
$f$ is a function of finite type, that $\cb \in \Sigma_f $ and that 
$\T_0$ is the polar type of $f$ at $\cb$.\\
In addition, we have 
$H_{\cb}(f;\s)=\left(\prod_{i=1}^n s_i^{u_i}\right) {\cal M}(f;\cb+\s)=
\prod_{i=1}^n \left(s_i^{u_i} Z(f_i;c_i+s_i)\right)$. Thus, 
$H_{\cb}(f;\zerob)=\prod_{i=1}^n A_i\neq 0$.
The corollary now follows by  theorem \ref{partieprincipale}.\cqfd
\subsection{Proof of corollary \ref{squarefree}}
It is easy to see that  corollary \ref{squarefree} follows from corollary \ref{mats+}.\\
We first observe that  $N(1_B;P;t)=N(\f;P;t)$ where 
each $f_j:\N\rightarrow \C$ is defined by 
$f_j(m)=1$ if $m$ is $k-$free and  $f_j(m)=0$ if not.\\
It is then a standard exercise to show that for all $j$ and for $\sigma >1$:
$Z(f_j;s)=\sum_{m=1}^{\infty}\frac{f_j(m)}{m^s}=\prod_p
(\sum_{\nu =0}^{k-1}\frac{1}{p^{\nu s}})=\frac{\zeta(s)}{\zeta(k s)}$. 
The assumptions of corollary \ref{mats+} are then satisfied with:\\
$\cb=(1,\dots,1)=\unb$, $I=\{\e_1,\dots ,\e_n\}$,  and 
$u_j=1,\ A_j=\frac{1}{\zeta(k)}$, \quad for all
$j=1,\dots, n$. \\ 
Moreover it is easy to see in our case that 
the mixed volume constant $A_0(\T_0;P)$ is equal to the Sargos 
constant $A_0(P)$. Thus,  corollary \ref{squarefree}
follows from corollary \ref{mats+} and the expression for $A_0(P) $  given by
  proposition \ref{sargoselliptique}. \cqfd
\subsubsection{Proof of corollary \ref{squarefreeuniform}}
We will derive  corollary \ref{squarefreeuniform} from corollary
\ref{applicationneuve}. Let 
$f:={\unb}_{D_k}$ be the characteristic function of $D_k$.
Set $I:=\{\e_1,\dots,\e_n\}$ and 
${\cal D}_f:=\{\s \in \C^n \mid \la \betab , \Re(\s)\ra >1 ~\forall
\betab \in I\}=\{\s \in \C^n \mid \sigma_i >1 ~\forall i\}$.\\
It is clear that ${\cal M}(f;\s):=\sum_{\m \in \N^{n}}
\frac{f(m_1,\dots,m_n)}{m^{s_1}\dots m_n^{s_n}}$ converges in 
${\cal D}_f$.  
The multiplicativity of $f$ implies also that for all
$\s \in {\cal D}_f$: 
$\ds {\cal M}(f;\s)=
\prod_p \left(\sum_{\{\nub \in \N_0^n; ~|\nub |\leq k-1\}} 
\frac{1}{p^{\la \nub ,\s\ra}}\right)$.\\
Let $\cb \in \partial {\cal D}_f=\{\s \in \C^n
\mid \sigma_i \geq 1 \forall i {\mbox { with at least one
    equality}}\}$.\\
Set $I_{\cb}:=\{\betab \in I\mid \la \betab, \cb\ra =1\}=\{\e_i \mid
c_i =1\}$ and $H_{\cb}(f;\s):=\left(\prod_{\betab \in I_{\cb}} \la
  \betab ,\s\ra \right) {\cal M}(f;\cb+\s)$.\\
It is easy to see that for all $\s \in \C^n$ satisfying $\sigma_i >0$
for all $i$, we have
\begin{equation}\label{Hcpositive}
H_{\cb}(f;\s)=\left(\prod_{\e_i \in I_{\cb}} s_i \zeta
  (1+s_i)\right) \prod_p \prod_{\e_i\in
  I_{\cb}}\left(1-\frac{1}{p^{1+s_i}}\right) \left(\sum_{|\nub |\leq k-1} 
\frac{1}{p^{\la \nub ,\cb\ra + \la \nub ,\s\ra}}\right)
\end{equation}
It follows from the standard properties of the Riemann zeta function and
the definition of $I_{\cb}$ that there exists $\eps_0>0$ such that $\s \mapsto
H_{\cb}(f;s)$ has a holomorphic continuation with moderate growth to
$\{\s \in \C^n \mid \sigma_i > -\eps_0 \forall i\}$ given by the right
side of the equality (\ref{Hcpositive}) above. Therefore $f$ is
a function of finite type.\par
Moreover ${\cal F}(\Sigma_f)(\unb)=\{\unb\}$ and $\cb=\unb \in
con^*\left(supp(P)\right)=\R_+^{*n}$. In addition, it is  easy to see that
assumption (1) and (2) of theorem \ref{partieprincipale} are
satisfied with $\cb =\unb$.\\
By analytic continuation, it follows also from (\ref{Hcpositive}) that 
$$H_{\unb}(f;\zerob)= \prod_p \left(1-\frac{1}{p}\right)^{n}\left(\sum_{|\nub |\leq k-1} 
\frac{1}{p^{|\nub |}}\right)= \prod_p \left(1-\frac{1}{p}\right)^{n}\left(1+\sum_{l=1}^{k-1} 
\frac{{n+l-1 \choose l}}{p^l}\right) >0.$$
We conclude that all three assumptions of theorem
\ref{partieprincipale} are satisfied. We can then apply corollary
\ref{applicationneuve} with:\\
$\cb=(1,\dots,1)=\unb$ and  $\T_{\cb} =(I_{\cb},\u),$ where $I_{\cb }
=\{\e_1,\dots ,\e_n\}$ and $u(\e_j)=1$ for all $j$.\\ 
Moreover it is easy to see in our case that $\iota(f)=|\unb|=n$, that 
$\rho (f;P)=\rho_0(\T_{\unb};P)=\rho_0(\T_{\unb})=1$ and that 
the mixed volume constant $A_0(\T_{\unb};P)$ is equal to the Sargos 
constant $A_0(P)$. Thus, corollary \ref{squarefreeuniform}
follows from corollary \ref{applicationneuve} and the expression for $A_0(P) $  given by
  proposition \ref{sargoselliptique}. \cqfd

\subsubsection{Proof of corollary \ref{goldba}}
Define for each $k=1,\dots,n$ and $i=1,2,$
$f_k^{(i)}: \N\rightarrow \R_+$
by $f_k^{(2)} =\Lambda$ and \\
$f_k^1(m)=\log (m)$ if $m$ is   prime  
and $f_k^{(1)}(m)=0$ if not.\\
It is well known that under $(RH)$ each of the
two Dirichlet  series
$s\mapsto \sum_{m\in \N} \frac{\Lambda(m)}{m^s}$ and
$s\mapsto \sum_{p\in {\cal P}} \frac{\log p}{p^s}$
has meromorphic continuation with moderate growth to the
half-plane $\{\Re (s)>\frac{1}{2}\}$ with a single  pole
at $s=1$ of order 1 and residue equal to 1.\\
Moreover, using the notations from corollary \ref{mats+},
we have $\T_0=(I;\u),$ where
$I=\{\e_1,\dots,\e_n\}$ and $u(\e_j)=1$ $\forall j$.
It follows that $A_0(\T_0;H)=A_0(H)$ is independent of $i$.   
Since $H$ is elliptic, the Sargos constant
$A_0(H)$ is given by proposition \ref{sargoselliptique}.
Thus corollary \ref{goldba} also follows from
corollary \ref{mats+}. \cqfd


\end{document}